\documentclass[11pt,reqno]{amsart}
\usepackage[a4paper,margin=30mm,headheight=14pt]{geometry}
\usepackage{amsmath,amssymb,amsthm,mathtools,mathrsfs,booktabs,array}
\usepackage{xcolor,enumitem,etoolbox,microtype,needspace}
\usepackage[T1]{fontenc}
\usepackage{lmodern}
\usepackage{hyperref}
\usepackage{aliascnt}
\usepackage[nameinlink,capitalize,noabbrev]{cleveref}
\hypersetup{
  colorlinks=true,
  linkcolor=blue!45!black,
  citecolor=blue!45!black,
  urlcolor=blue!45!black,
  pdftitle={Critical Nonlinear Potential Inequalities: Path Testing and Anisotropic Cancellation},
  pdfauthor={Zhouyu Long and Wenming Zou},
  pdfsubject={Spatial nonlinearities, vector mass paths, and anisotropic critical potential inequalities},
  pdfkeywords={critical nonlinear potential inequalities, homogeneous kernels, path testing, vector measures, spatial nonlinearities, measurable targets, anisotropic dilations, Maz'ya--Phi inequalities, cross-scale interaction}
}
\numberwithin{equation}{section}

\newtheorem{theorem}{Theorem}[section]
\theoremstyle{plain}
\newaliascnt{lemma}{theorem}
\newtheorem{lemma}[lemma]{Lemma}
\aliascntresetthe{lemma}
\theoremstyle{plain}
\newaliascnt{proposition}{theorem}
\newtheorem{proposition}[proposition]{Proposition}
\aliascntresetthe{proposition}
\theoremstyle{plain}
\newaliascnt{corollary}{theorem}
\newtheorem{corollary}[corollary]{Corollary}
\aliascntresetthe{corollary}
\theoremstyle{definition}
\newaliascnt{definition}{theorem}

\aliascntresetthe{definition}
\theoremstyle{remark}
\newaliascnt{remark}{theorem}
\newtheorem{remark}[remark]{Remark}
\aliascntresetthe{remark}
\theoremstyle{remark}
\newaliascnt{example}{theorem}
\newtheorem{example}[example]{Example}
\aliascntresetthe{example}
\theoremstyle{plain}
\crefname{theorem}{Theorem}{Theorems}
\Crefname{theorem}{Theorem}{Theorems}
\crefname{lemma}{Lemma}{Lemmas}
\Crefname{lemma}{Lemma}{Lemmas}
\crefname{proposition}{Proposition}{Propositions}
\Crefname{proposition}{Proposition}{Propositions}
\crefname{corollary}{Corollary}{Corollaries}
\Crefname{corollary}{Corollary}{Corollaries}
\crefname{definition}{Definition}{Definitions}
\Crefname{definition}{Definition}{Definitions}
\crefname{remark}{Remark}{Remarks}
\Crefname{remark}{Remark}{Remarks}
\crefname{example}{Example}{Examples}
\Crefname{example}{Example}{Examples}
\setlist[enumerate,1]{label=\textup{(\roman*)},ref=\textup{(\roman*)}}

\newcommand{\R}{\mathbb R}
\newcommand{\supp}{\operatorname{supp}}
\newcommand{\dist}{\operatorname{dist}}
\newcommand{\dd}{\,\mathrm d}
\newcommand{\one}{\mathbf 1}

\title[Critical nonlinear potential inequalities]
{Critical Nonlinear Potential Inequalities:\\
Path Testing and Anisotropic Cancellation}
\author[Zhouyu Long]{\bfseries Zhouyu Long$^{*}$}
\author[Wenming Zou]{\bfseries Wenming Zou$^{\dagger}$}
\thanks{$^{*}$Department of Mathematical Sciences, Tsinghua University,
Beijing 100084, China.\newline
E-mail: \href{mailto:longzy25@mails.tsinghua.edu.cn}{\nolinkurl{longzy25@mails.tsinghua.edu.cn}}.}
\thanks{$^{\dagger}$Department of Mathematical Sciences, Tsinghua University,
Beijing 100084, China.\newline
E-mail: \href{mailto:zou-wm@mail.tsinghua.edu.cn}{\nolinkurl{zou-wm@mail.tsinghua.edu.cn}}.}

\makeatletter
\patchcmd{\@setauthors}{\centering\footnotesize}{\centering\large}{}%
  {\ClassError{amsart}{Unable to set author font}{Check the amsart author layout.}}
\patchcmd{\@setauthors}{\MakeUppercase{\authors}}{\authors}{}%
  {\ClassError{amsart}{Unable to preserve author case}{Check the amsart author layout.}}
\def\author@andify{\nxandlist{\unskip,\quad}{\unskip,\quad}{\unskip,\quad}}
\let\addresses\@empty
\makeatother
\date{}

\begin{document}
\begin{abstract}
We characterize critical Maz'ya--$\Phi$ inequalities for compactly supported
vector measures with locally $L^p$ potentials. The nonlinearity may depend
measurably on position and is positively $p$-homogeneous and uniformly
Lipschitz on the sphere. For geometrically shrinking nested symmetric convex
shells with critical kernel size and a Dini $L^1$ translation modulus, the
mean-zero estimate is equivalent to a uniform bound on finite vector mass
paths. For scalar sources and bounded target weights, the test reduces to
signed point-mass annular sums. Arbitrary total mass also requires a finite
absolute coarse-scale tail, automatic on finite-measure targets. Applications
include exact criteria for measurable targets and higher-order
constant-coefficient elliptic estimates, a proof of Stolyarov's Conjecture~5.1
for fixed positive diagonal dilations, and its extension to expanding matrix
dilations using Lyapunov sections.
\end{abstract}
\maketitle

\section{Introduction and main results}\label{sec:introduction}

\subsection{Background and the critical point-mass obstruction}

Let $d\ge 2$, $0<\alpha<d$, and let $p=d/(d-\alpha)$ be the critical
Hardy--Littlewood--Sobolev exponent. For a kernel
$K(r\omega)=r^{\alpha-d}\widetilde K(\omega)$ on $\R^d\setminus\{0\}$ taking
values in $\R^\ell$ and a positively $p$-homogeneous nonlinearity
$\Phi:\R^\ell\to\R$, substituting a signed Dirac mass $\sigma\delta_0$
($\sigma=\pm1$) yields the scale-invariant radial density
\[
 \Phi(\sigma K(r\omega))\,r^{d-1}\dd r
 =\Phi(\sigma\widetilde K(\omega))\,\frac{\dd r}{r}.
\]
Because the Haar measure $\dd r/r$ has infinite mass both near the origin and
at infinity, signed accumulation across scales constitutes the fundamental
obstruction to critical $L^1$ potential inequalities.

Motivated by this critical mechanism, Maz'ya~\cite[Section~2]{Mazya} posed the
gradient--Laplacian problem on $\R^d$, seeking cancellation conditions on
$\Phi$ for which
\[
 \left|\int_{\R^d}\Phi(\nabla u)\,\dd x\right|
 \le C\|\Delta u\|_{L^1(\R^d)}^{d/(d-1)}
 \qquad\bigl(u\in C_c^\infty(\R^d;\R)\bigr).
\]
Guided by a martingale model~\cite{StolyarovMartingale},
Stolyarov~\cite[Theorem~2.1]{StolyarovWhole} resolved Maz'ya's conjecture and
characterized the general Euclidean whole-space inequality in terms of the two
full-sphere cancellation conditions
$\int_{\mathbb S^{d-1}}\Phi(\pm\widetilde K(\omega))\,\dd\sigma(\omega)=0$.
On bounded $C^{1,\beta_\partial}$ domains ($0<\beta_\partial\le 1$),
Stolyarov~\cite[Theorems~1.2--1.3]{StolyarovDomains} subsequently showed that
cancellation over every tangent hemisphere is both necessary and sufficient
(and automatically implies full-sphere cancellation), whereas on a fixed
half-space both hemisphere and full-sphere cancellation are required.

We give a necessary and sufficient test on arbitrary measurable targets,
extending the signed planar quadratic criterion of Long and
Zou~\cite[Theorem~1.2]{LongZouCumulative}. The same argument proves the
anisotropic whole-space classification conjectured by
Stolyarov~\cite[Conjecture~5.1]{StolyarovWhole}. Both follow from a path
classification for spatial nonlinearities on nested convex shells.

\subsection{Measurable targets and anisotropic cancellation}
\label{sec:intro-headline}

We use Euclidean norms on finite-dimensional spaces and the induced operator
norm on spaces of linear maps. Let $d\ge 2$, $\ell\ge 1$, and $p>1$. For a
positively $p$-homogeneous function
$\Phi:\R^\ell\to\R$ with $\Phi(0)=0$, we measure its regularity on the unit
sphere $\mathbb S^{\ell-1}$ (with the induced Euclidean metric, including
$\mathbb S^0=\{-1,1\}$ when $\ell=1$) by
\[
 L_\Phi:=\|\Phi\|_{L^\infty(\mathbb S^{\ell-1})}
 +\operatorname{Lip}_{\mathbb S^{\ell-1}}(\Phi)<\infty.
\]
For positively $p$-homogeneous functions, finiteness of $L_\Phi$ is equivalent
to local Lipschitz continuity on $\R^\ell$. Optimal constants in all
inequalities below are understood as infima over admissible nonnegative
constants, with $\inf\varnothing=+\infty$.

We first treat the isotropic Euclidean setting. Let $\alpha\in(0,d)$,
$p=d/(d-\alpha)$, and $K(x)=|x|^{\alpha-d}\widetilde K(x/|x|)$ for $x\ne 0$,
with $K(0)=0$ and $\widetilde K:\mathbb S^{d-1}\to\R^\ell$ Lipschitz. For
every $f\in L_c^\infty(\R^d;\R)$, the potential $K*f$ belongs to
$C_0(\R^d;\R^\ell)$ (Section~\ref{sec:shell-geometry}), so $\Phi(K*f)$ is
bounded and integrable on every measurable set $E\subset\R^d$ of finite
measure. To capture signed cross-scale accumulation on $E$, we introduce the
continuous signed point-mass seminorm
\begin{equation}\label{eq:continuous-pointmass-cocycle}
 \mathfrak P_E^{\rm cont}
 :=\sup_{\substack{a\in\R^d,\ \sigma=\pm 1\\0<r<R<\infty}}
 \left|\int_{E\cap\{r<|x-a|<R\}}\Phi(\sigma K(x-a))\,\dd x\right|,
\end{equation}
which depends only on $(K,\Phi)$ and the Lebesgue equivalence class of $E$.
Let $C_E^*$ denote the optimal constant in
\begin{equation}\label{eq:sharp-measurable-density}
 \left|\int_E\Phi(K*f)\,\dd x\right|\le C\|f\|_1^p
 \qquad\bigl(f\in L_c^\infty(\R^d;\R)\bigr).
\end{equation}

\Needspace{9\baselineskip}
\begin{theorem}
\label{thm:sharp-measurable-characterization}
Let $E\subset\R^d$ be a measurable set with $0<|E|<\infty$. Then
$C_E^*<\infty$ if and only if $\mathfrak P_E^{\rm cont}<\infty$. Moreover,
\begin{equation}\label{eq:main-continuous-quantitative}
 C^{-1}\bigl(1+\mathfrak P_E^{\rm cont}\bigr)\le 1+C_E^*
 \le C\bigl(1+\mathfrak P_E^{\rm cont}\bigr),
\end{equation}
where $C\ge 1$ depends only on $d$, $\alpha$, $\ell$,
$\|\widetilde K\|_\infty+\operatorname{Lip}(\widetilde K)$, and $L_\Phi$,
independently of $E$.
\end{theorem}

An equivalent dyadic annular test is established in
Lemma~\ref{lem:continuous-dyadic-cocycle}. On arbitrary measurable sets $E$
(including $|E|=\infty$), mean-zero sources remain governed solely by
$\mathfrak P_E^{\rm cont}<\infty$, while unconstrained sources require in
addition a finite large-scale absolute point-mass tail
(Corollary~\ref{cor:weighted-homogeneous-target}).

We next turn to anisotropic dilations. Fix a positive weight vector
$\mathbf a=(a_1,\ldots,a_d)\in(0,\infty)^d$ normalized by $\sum_{j=1}^d a_j=d$,
and set $\delta_t^{\mathbf a}x=(t^{a_j}x_j)_{j=1}^d$ for $t>0$. Define the
associated anisotropic radius $\rho_{\mathbf a}:\R^d\to[0,\infty)$, smooth
off the origin, as the unique positive solution for $x\ne 0$ of
\begin{equation}\label{eq:aniso-radius}
 \sum_{j=1}^d x_j^2\,\rho_{\mathbf a}(x)^{-2a_j}=1,
 \qquad \rho_{\mathbf a}(0)=0,
\end{equation}
so that $\rho_{\mathbf a}(\delta_t^{\mathbf a}x)=t\rho_{\mathbf a}(x)$.
Set $\beta=d/p$ and $\alpha=d-\beta\in(0,d)$. Let
$\mathcal K:\R^d\setminus\{0\}\to\R^\ell$ satisfy
$\mathcal K(\delta_t^{\mathbf a}x)=t^{-\beta}\mathcal K(x)$ for $t>0$, with
$\mathcal K(0)=0$ and Lipschitz spherical restriction
$\widetilde{\mathcal K}=\mathcal K|_{\mathbb S^{d-1}}$. Let
$w_{\mathbf a}(\omega):=\sum_{j=1}^d a_j\omega_j^2$ denote the anisotropic
polar flux weight (Lemma~\ref{lem:aniso-polar}). Every mean-zero
$f\in C_c^\infty(\R^d;\R)$ satisfies $\mathcal K*f\in L^p(\R^d;\R^\ell)$, so
the integral below converges absolutely.

\Needspace{9\baselineskip}
\begin{theorem}
\label{thm:anisotropic-mazya-phi}
Under the preceding hypotheses, the following statements are equivalent:
\begin{enumerate}
\item There exists a finite constant $C$ such that
\begin{equation}\label{eq:anisotropic-mazya-phi}
 \left|\int_{\R^d}\Phi(\mathcal K*f)\,\dd x\right|
 \le C\|f\|_1^p
 \qquad\Bigl(f\in C_c^\infty(\R^d;\R),\ \int_{\R^d} f=0\Bigr).
\end{equation}
\item For both signs $\sigma=\pm 1$,
\begin{equation}\label{eq:anisotropic-mazya-cancellations}
 \int_{\mathbb S^{d-1}}w_{\mathbf a}(\omega)\,
 \Phi\bigl(\sigma\widetilde{\mathcal K}(\omega)\bigr)\,\dd\sigma(\omega)=0.
\end{equation}
\end{enumerate}
The admissible constant $C$ depends only on $d$, $\ell$, $\mathbf a$, $\alpha$,
and the spherical supremum norms and Lipschitz seminorms of
$\widetilde{\mathcal K}$ and $\Phi$.
\end{theorem}

Theorem~\ref{thm:anisotropic-mazya-phi} proves Stolyarov's
Conjecture~5.1~\cite[(5.5)--(5.7)]{StolyarovWhole} for every fixed positive
diagonal weight vector $\mathbf a$ in $d\ge 2$; the estimate is uniform over
mean-zero sources, with a constant that may depend on $\mathbf a$ (and may
deteriorate as some $a_j\downarrow 0$). In dimension $d=1$, the normalization
forces $a_1=1$, which is the isotropic case of
\cite[Theorem~2.1]{StolyarovWhole}. Generalizations to non-diagonal expanding
real matrices $T_t=t^A$ in Lyapunov sections and to weighted measurable targets
are established in Section~\ref{sec:anisotropic-measurable}.

\subsection{The unified framework: convex shells, spatial nonlinearities, and path testing}
\label{sec:intro-general}

Both Theorem~\ref{thm:sharp-measurable-characterization} and
Theorem~\ref{thm:anisotropic-mazya-phi} emerge from a general classification
for vector-valued measure sources, spatially dependent nonlinearities, and
nonstationary convex-shell kernel families.

Let $(A_n)_{n\in\mathbb Z}$ be a nested sequence of closed, origin-symmetric
convex bodies in $\R^d$ with nonempty interiors satisfying the uniform
geometric contraction
\begin{equation}\label{eq:common-scale-hypotheses}
 A_{n+1}\subset A_n,\qquad
 A_{n+k}\subset C_s\,2^{-b_*k}A_n\qquad (k\ge 0),
\end{equation}
for fixed constants $C_s\ge 1$ and $b_*>0$. Write $V_n:=|A_n|$, let
$D_n(h):=\inf\{t>0:h\in tA_n\}$ be the Minkowski norm of $A_n$, and define the
disjoint shells $E_n:=\operatorname{int}A_n\setminus\operatorname{int}A_{n+1}$.
Let $\vartheta:[0,\infty)\to[0,\infty)$ be a bounded, nondecreasing, continuous
modulus with $\vartheta(0)=0$ satisfying the Dini condition
\begin{equation}\label{eq:dini-modulus}
 \mathcal D_\vartheta:=\|\vartheta\|_\infty
 +\int_0^1\vartheta(t)\,\frac{\dd t}{t}<\infty.
\end{equation}
Fix $m\ge 1$. We consider measurable operator-valued kernels
$H_n:\R^d\to\mathcal L(\R^m,\R^\ell)$ ($n\in\mathbb Z$) satisfying the
critical size bound and Dini $L^1$ translation modulus
\begin{equation}\label{eq:common-shell-hypotheses}
 \begin{gathered}
  \{H_n\ne 0\}\subset E_n,\qquad
  \|H_n\|_\infty\le C_0 V_n^{-1/p},\\
  \|H_n(\cdot-h)-H_n\|_1\le C_1 V_n^{1/p'}\vartheta\bigl(D_n(h)\bigr)
  \qquad (n\in\mathbb Z,\ h\in\R^d),
 \end{gathered}
\end{equation}
where $p'=p/(p-1)$ and all norms use the operator norm on
$\mathcal L(\R^m,\R^\ell)$. Choosing bounded Borel representatives that
enforce the support and size bounds pointwise, we set
$K_H:=\sum_{n\in\mathbb Z} H_n$. In the diagonal dilation setting, the
ellipsoids $A_n=\delta_{2^{-n}}^{\mathbf a}\overline{B(0,1)}$ give
$V_n=|B(0,1)|\,2^{-dn}$, $D_n(h)=|\delta_{2^n}^{\mathbf a}h|$, and
\begin{equation}\label{eq:annular-family-hypotheses}
 \begin{aligned}
  \{H_n\ne 0\}&\subset\{2^{-n-1}\le\rho_{\mathbf a}<2^{-n}\},\\
  \|H_n\|_\infty&\le C_0\,2^{\beta n},\\
  \|H_n(\cdot-h)-H_n\|_1&\le C_1\,2^{-\alpha n}\vartheta\bigl(|\delta_{2^n}^{\mathbf a}h|\bigr).
 \end{aligned}
\end{equation}

Next, let $\Psi:\R^d\times\R^\ell\to\R$ be measurable in $x$ for each $z$.
Assume that on a common full-measure set of $x\in\R^d$, $\Psi(x,0)=0$,
$\Psi(x,tz)=t^p\Psi(x,z)$ for all $t>0$, and
\begin{equation}\label{eq:spatial-nonlinearity-norm}
 L_\Psi:=\operatorname*{ess\,sup}_{x\in\R^d}
 \Bigl(\|\Psi(x,\cdot)\|_{L^\infty(\mathbb S^{\ell-1})}
 +\operatorname{Lip}_{\mathbb S^{\ell-1}}\Psi(x,\cdot)\Bigr)<\infty.
\end{equation}
Selecting Borel versions on a countable dense subset of $\R^\ell$, extending
by continuity off one common Borel null set of $x$, and setting $\Psi=0$ on
that null set produces a jointly Borel Carath\'eodory representative of $\Psi$.

For a finite $\R^m$-valued Borel measure $\nu$, let $|\nu|$ denote its total
variation measure with respect to the Euclidean norm on $\R^m$, and set
$\|\nu\|_{\rm TV}:=|\nu|(\R^d)$. Let $\mathcal X_H$ denote the class of
compactly supported finite $\R^m$-valued Borel measures $\nu$ satisfying
$K_H*\nu\in L^p_{\mathrm{loc}}(\R^d;\R^\ell)$. This class contains
$L_c^\infty(\R^d;\R^m)$ and compactly supported singular measures with finite
small-scale energy (Section~\ref{sec:concentrated-measures}).

For $a\in\R^d$, $n\in\mathbb Z$, and $v\in\R^m$, define the single-shell
point-mass response
\[
 q_{n,a}(v):=\int_{\R^d}\Psi\bigl(x,H_n(x-a)v\bigr)\,\dd x.
\]
Given a finite path $v_N,\ldots,v_M\in\R^m$ ($N\le M$), its zero-extended total
variation is
\[
 D(v_N,\ldots,v_M)
 :=|v_N|+\sum_{n=N}^{M-1}|v_{n+1}-v_n|+|v_M|.
\]
Define the path-testing constant and large-scale directional tail by
\begin{align}
 \mathfrak V_\Psi^H
 &:=\sup_{\substack{a\in\R^d,\ N\le M\\v_N,\ldots,v_M\in\R^m,\ D>0}}
 \frac{\bigl|\sum_{n=N}^M q_{n,a}(v_n)\bigr|}{D(v_N,\ldots,v_M)^p},
 \label{eq:vector-path-primitive}\\
 \mathfrak T_\Psi^H(a_0,v)
 &:=\sum_{n\le 0}\int_{\R^d}\bigl|\Psi\bigl(x,H_n(x-a_0)v\bigr)\bigr|\,\dd x.
 \label{eq:vector-directional-tail}
\end{align}

\begin{theorem}
\label{thm:path-family-classification}
Assume \eqref{eq:common-scale-hypotheses}, \eqref{eq:dini-modulus},
\eqref{eq:common-shell-hypotheses}, and \eqref{eq:spatial-nonlinearity-norm}.
\begin{enumerate}
\item For every $\nu\in\mathcal X_H$ with $\nu(\R^d)=0$, the integrand
$\Psi(\cdot,K_H*\nu)$ belongs to $L^1(\R^d)$. The mean-zero estimate
\begin{equation}\label{eq:vector-zero-mean-estimate}
 \left|\int_{\R^d}\Psi\bigl(x,K_H*\nu(x)\bigr)\,\dd x\right|
 \le C_{\Psi,0}\|\nu\|_{\rm TV}^p
 \qquad\bigl(\nu\in\mathcal X_H,\ \nu(\R^d)=0\bigr)
\end{equation}
holds with a finite constant if and only if $\mathfrak V_\Psi^H<\infty$, and
its optimal constant $C_{\Psi,0}^{H,*}$ satisfies
\begin{equation}\label{eq:vector-zero-mean-quantitative}
 L_\Psi+C_{\Psi,0}^{H,*}\asymp L_\Psi+\mathfrak V_\Psi^H.
\end{equation}
\item Absolute integrability of $\Psi(\cdot,K_H*\nu)$ for every
$\nu\in\mathcal X_H$ is equivalent to
\begin{equation}\label{eq:vector-all-direction-tails}
 \mathfrak T_\Psi^H(a_0,v)<\infty
 \qquad\text{for every }v\in\R^m.
\end{equation}
Under \eqref{eq:vector-all-direction-tails}, the bound in
\eqref{eq:vector-zero-mean-estimate} holds for all $\nu\in\mathcal X_H$,
regardless of total mass, if and only if $\mathfrak V_\Psi^H<\infty$. Its
optimal constant $C_\Psi^{H,*}$ satisfies
\begin{equation}\label{eq:vector-full-quantitative}
 L_\Psi+C_\Psi^{H,*}\asymp L_\Psi+\mathfrak V_\Psi^H.
\end{equation}
\end{enumerate}
In each integrable branch, restricting the source class to
$C_c^\infty(\R^d;\R^m)$ preserves the optimal constant. For each $v\in\R^m$, finiteness of
$\mathfrak T_\Psi^H(a_0,v)$ is independent of the center $a_0$ and of any
fixed scale cutoff. The implicit comparison constants depend only on the
structural data $(d,m,\ell,p,C_s,b_*,C_0,C_1,\mathcal D_\vartheta)$,
independently of $\Psi$ and of finite tail values.
\end{theorem}

By Lemma~\ref{lem:scalar-path-reduction}, finiteness of
\eqref{eq:vector-path-primitive} is equivalent to uniform testing on
unit-vector paths of interior variation at most one. Fixed-direction testing
can fail for vector-valued sources even when the kernels have smooth angular
profiles and the nonlinearity has a smooth restriction to the sphere
(Proposition~\ref{prop:directional-testing-fails}).

For scalar sources ($m=1$), we identify $H_n$ with a map into $\R^\ell$ and
set $\Psi(x,z)=b(x)\Phi(z)$, where $b\in L^\infty(\R^d;\R)$ with
$B:=\|b\|_\infty$. Introducing
\begin{align}
 \mathfrak P_b^H
 &:=\sup_{\substack{a\in\R^d,\ \sigma=\pm 1\\N\le M}}
 \left|\sum_{n=N}^M\int_{\R^d}b(x)\Phi\bigl(\sigma H_n(x-a)\bigr)\,\dd x\right|,
 \label{eq:family-point-mass-primitive}\\
 \Lambda_H
 &:=\sup_{\sigma=\pm 1,\ N\le M}
 \left|\sum_{n=N}^M\int_{\R^d}\Phi\bigl(\sigma H_n(x)\bigr)\,\dd x\right|,
 \label{eq:family-whole-primitive}\\
 \mathfrak T_b^H(a_0)
 &:=\max_{\sigma=\pm 1}\sum_{n\le 0}\int_{\R^d}
 |b(x)|\,\bigl|\Phi\bigl(\sigma H_n(x-a_0)\bigr)\bigr|\,\dd x,
 \label{eq:family-absolute-tail}
\end{align}
Lemma~\ref{lem:scalar-path-reduction} yields $\mathfrak V_{b\Phi}^H=2^{-p}\mathfrak P_b^H$
and produces the following scalar weighted classification.

\begin{corollary}
\label{thm:annular-family-classification}
Assume the kernel hypotheses of Theorem~\ref{thm:path-family-classification}
with $m=1$, and let $b\in L^\infty(\R^d;\R)$ with $B=\|b\|_\infty$.
\begin{enumerate}
\item For every mean-zero $f\in L_c^\infty(\R^d;\R)$, the integral in
\begin{equation}\label{eq:family-zero-mean-estimate}
 \left|\int_{\R^d}b(x)\Phi(K_H*f(x))\,\dd x\right|\le C_{b,0}\|f\|_1^p
\end{equation}
converges absolutely. The estimate holds uniformly on this class if and only
if $\mathfrak P_b^H<\infty$, with optimal constant $C_{b,0}^{H,*}$ satisfying
\begin{equation}\label{eq:family-zero-mean-quantitative}
 B+C_{b,0}^{H,*}\asymp B+\mathfrak P_b^H.
\end{equation}
\item Absolute convergence for every $f\in L_c^\infty(\R^d;\R)$ is equivalent
to $\mathfrak T_b^H(a_0)<\infty$. Under this tail condition, the unconstrained
estimate
\begin{equation}\label{eq:family-finite-target-estimate}
 \left|\int_{\R^d}b(x)\Phi(K_H*f(x))\,\dd x\right|\le C_b\|f\|_1^p
 \qquad\bigl(f\in L_c^\infty(\R^d;\R)\bigr)
\end{equation}
holds uniformly if and only if $\mathfrak P_b^H<\infty$, with optimal constant
$C_b^{H,*}$ satisfying
\begin{equation}\label{eq:family-finite-quantitative}
 B+C_b^{H,*}\asymp B+\mathfrak P_b^H.
\end{equation}
The tail condition holds automatically whenever $b\in L^1(\R^d)\cap L^\infty(\R^d)$.
\end{enumerate}
When $b\equiv 1$, translation invariance gives $\mathfrak P_1^H=\Lambda_H$, so
the whole-space estimate
\begin{equation}\label{eq:family-whole-space-estimate}
 \left|\int_{\R^d}\Phi(K_H*f(x))\,\dd x\right|\le C_{\rm wh}\|f\|_1^p
 \qquad\Bigl(f\in C_c^\infty(\R^d;\R),\ \int_{\R^d}f=0\Bigr)
\end{equation}
holds if and only if $\Lambda_H<\infty$, with
$1+C_{\rm wh}^{H,*}\asymp 1+\Lambda_H$. Both optimal constants in
\textup{(i)} and \textup{(ii)} are unchanged on $C_c^\infty(\R^d;\R)$. If $B$
is replaced by $BL_\Phi$ on both sides, the comparison constants depend only on
the structural data of Theorem~\ref{thm:path-family-classification}; the
displayed comparisons may additionally depend on $\max(1,L_\Phi)$,
independently of $b$, the scale endpoints, and finite tail values.
\end{corollary}

Theorem~\ref{thm:path-family-classification} and
Corollary~\ref{thm:annular-family-classification} also hold in dimension $d=1$:
the slicing step in \eqref{eq:shell-collar} uses the interval identity
$|I\mathbin\triangle(I+h)|=2\min\{|h|,|I|\}$. The transport and path arguments
then apply in dimension one.

\subsection{Proof strategy and organization of the paper}

The Euclidean proofs in \cite{StolyarovWhole,StolyarovDomains} use dyadic
cube-mass energies and the Three Lattice Theorem, guided by the martingale
model~\cite{StolyarovMartingale}. Transferring this localization to anisotropic
dilations is one of the difficulties identified in
\cite[Section~5]{StolyarovWhole}. Here directional slicing of convex bodies
gives the shell collar estimate \eqref{eq:shell-collar}. Source-centered
normalization and transport then control both the error of replacing each
shell response by a point-mass response and the nonlinear interaction between
scales. Concentric uniform densities realize finite mass paths, and
large-scale point-mass comparison controls convergence and mollification.

Sections~\ref{sec:pointmass-energy}--\ref{sec:proof-sharp} prove the convex-shell
classification. Sections~\ref{sec:anisotropic-measurable} and
\ref{sec:pde-and-extensions} treat homogeneous dilations, measurable targets,
and elliptic estimates. Section~\ref{sec:sharpness} gives counterexamples to
several relaxations and separates signed block tests from absolute shell sums.

\subsection{Relation to previous work}

For bounded planar open sets $\Omega$, Long and
Zou~\cite[Theorem~1.1]{LongZouCumulative} prove
$\sum_n|\int_\Omega\Phi(K_n*f)\,\dd x|\lesssim\|f\|_1^p$ for $f\in L_c^1$
under plane cancellation and absolute summability of boundary-centered shell
responses, uniformly in the boundary center. For the Newton kernel and
trace-free quadratic forms, their Theorem~1.2 gives a signed all-center
characterization with an additive comparison of best constants.
Their Corollary~1.10 gives a tangent-cone criterion on finitely cornered Dini
domains for general admissible pairs.

On the common source class $L_c^\infty$, Theorem~\ref{thm:sharp-measurable-characterization}
extends the signed characterization to finite-measure targets in every
dimension $d\ge 2$ and general critical homogeneous kernels and nonlinearities
with Lipschitz angular parts. In the planar quadratic case on bounded targets,
an annulus is a difference of exterior truncations, so the signed tests agree
up to a factor of two. The planar PDE results in
\cite[Theorems~1.5 and~1.7]{LongZouCumulative} use
$\|\Delta u\|_{L^1(\Omega)}+\|\partial_n u\|_{L^1(\partial\Omega)}$.
Corollaries~\ref{cor:elliptic-measurable-target}
and~\ref{cor:laplacian-measurable-target} give criteria with
$\|Pu\|_{L^1(\R^d)}$ for compactly supported smooth functions.

Linear antecedents include Bourgain--Brezis Hodge
estimates~\cite{BourgainBrezis,BourgainBrezisFull}, Van Schaftingen's
ellipticity and canceling criterion~\cite[Theorem~1.3]{VanSchaftingen}, and
Stolyarov's Dirac-exclusion theorem on closed translation- and
dilation-invariant subspaces~\cite[Theorem~1]{StolyarovEndpoint}; see also the
surveys and domain or trace results
\cite{CianchiMazyaArbitraryDomains,GmeinederRaitaVanSchaftingen,RaitaSpectorStolyarovTrace,Spector,VanSchaftingenSurvey}.
Homogeneous dilations and stratified-group endpoint estimates appear
in~\cite{FollandSteinHomogeneousGroups,VanSchaftingenYung}. Recently,
Stolyarov~\cite[Theorem~1.3]{StolyarovAnisotropicBB} obtained linear
anisotropic Riesz potential estimates $L^1\to L^{d/(d-\alpha)}$ for $L^1$
fields in closed translation- and anisotropic-dilation-invariant subspaces
excluding vector Dirac masses. Our signed classification tests point masses
shell by shell to control compactly supported sources with locally $L^p$
potentials.

\section{Shell geometry and centered source transport}
\label{sec:pointmass-energy}

In this section we record the analytic properties of convex-shell kernel
families, derive the nonlinear increment and single-shell point-mass
replacement estimates, and prove that finite sums of single-shell integrals
are bounded by the path-testing constant $\mathfrak V_\Psi^H$.

\subsection{Convex-shell geometry and convolution properties}
\label{sec:shell-geometry}

Let $(A_n)_{n\in\mathbb Z}$ be the nested origin-symmetric convex bodies,
volumes $V_n=|A_n|$, Minkowski norms $D_n$, shells
$E_n=\operatorname{int}A_n\setminus\operatorname{int}A_{n+1}$, and kernels
$H_n$ satisfying \eqref{eq:common-scale-hypotheses}--\eqref{eq:common-shell-hypotheses}.
Because each $A_n$ is a closed, origin-symmetric convex body with nonempty
interior, each $D_n$ is a norm on $\R^d$, and geometric contraction ensures
that $(E_n)_{n\in\mathbb Z}$ partitions $\R^d\setminus\{0\}$. Since convex
boundaries have zero Lebesgue measure, changing shell boundary conventions
leaves convolutions unchanged almost everywhere by Tonelli's theorem.

By \eqref{eq:common-scale-hypotheses}, the volumes contract geometrically:
\[
 \frac{V_{n+k}}{V_n}\le C_s^d\,2^{-db_*k}\qquad(k\ge 0).
\]
Together with the support and $L^\infty$ bounds in
\eqref{eq:common-shell-hypotheses}, this yields
\begin{equation}\label{eq:annular-family-size}
 \|H_n\|_1\lesssim V_n^{1/p'},\qquad \|H_n\|_p\lesssim 1.
\end{equation}
Put $H_{\le n}:=\sum_{j\le n}H_j$ and $K_H:=\sum_{j\in\mathbb Z}H_j$. The
large-scale series $\sum_{j\le n}H_j$ converges uniformly (since $V_j^{-1/p}$
decays geometrically as $j\to-\infty$), while the small-scale series
$\sum_{j>n}H_j$ converges in $L^1(\R^d)$ (since $V_j^{1/p'}$ decays
geometrically as $j\to+\infty$). Hence $K_H$ is locally integrable and bounded
away from the origin.

For any finite $\R^m$-valued Borel measure $\nu$, the convolution $K_H*\nu$ is
well-defined almost everywhere by Tonelli's theorem. Moreover, for any
$\chi\in C_c^\infty(\R^d)$ and any bounded measurable set $S\subset\R^d$,
splitting $K_H=H_{\le 0}+(K_H-H_{\le 0})$ gives
\[
 \begin{aligned}
 &\int_S\int_{\R^d}\int_{\R^d}|\chi(h)|\,\|K_H(x-y-h)\|\,
              \dd|\nu|(y)\,\dd h\,\dd x\\
 &\qquad\le \|\chi\|_1\|\nu\|_{\rm TV}
    \bigl(|S|\|H_{\le 0}\|_\infty+\|K_H-H_{\le 0}\|_1\bigr)<\infty,
 \end{aligned}
\]
so Fubini's theorem yields the mollification associativity identity
\begin{equation}\label{eq:mollification-associativity}
 K_H*(\nu*\chi)=(K_H*\nu)*\chi\qquad\text{almost everywhere on }\R^d.
\end{equation}
For $f\in L_c^\infty(\R^d;\R^m)$, the large-scale $L^\infty$ and small-scale
$L^1$ bounds imply uniform convergence of finite-block convolutions to
$K_H*f$; since each $H_n*f$ is continuous with compact support,
$K_H*f\in C_0(\R^d;\R^\ell)$.

Throughout Sections~\ref{sec:pointmass-energy}--\ref{sec:proof-sharp}, the
structural data are $(d,m,\ell,p,C_s,b_*,C_0,C_1,\mathcal D_\vartheta)$, and
implicit constants depend only on these parameters unless stated otherwise;
dependence on the nonlinearity is displayed explicitly through $L_\Psi$. We
repeatedly use the geometric Dini summation bound
\begin{equation}\label{eq:dini-geometric-sum}
 \sum_{j\ge 0}\vartheta\bigl(C\,2^{-bj}\bigr)\lesssim_{C,b}\mathcal D_\vartheta
 \qquad(C,b>0),
\end{equation}
obtained by comparing each term with the integral of $\vartheta(t)\,\dd t/t$
over the preceding geometric interval.

\subsection{Nonlinear increments and single-shell point-mass replacement}

For $s,t\ge 0$, define the convexity increment
$\Delta_p(s,t):=(s+t)^p-s^p-t^p$. It is symmetric, nonnegative, nondecreasing
in each variable, and satisfies
\[
 \Delta_p(s,t)\asymp_p \min(s,t)\max(s,t)^{p-1}.
\]
Indeed, when $0\le t\le s$, convexity of $u\mapsto u^p$ at $u=s$ gives
$(s+t)^p-s^p\ge ps^{p-1}t$, and $t^p\le s^{p-1}t$ yields the lower bound
$(p-1)s^{p-1}t$; the upper bound follows from the mean-value theorem. In
particular,
\[
 \Delta_p(s,t)\lesssim_p\min(s^{p-1}t,st^{p-1})\quad(1<p\le 2),
 \qquad
 \Delta_p(s,t)\lesssim_p s^{p-1}t+st^{p-1}\quad(p\ge 2).
\]

By positive $p$-homogeneity and the spherical Lipschitz bound
\eqref{eq:spatial-nonlinearity-norm}, the integrand $\Psi$ satisfies
\begin{align}
 |\Psi(x,z)|&\le L_\Psi|z|^p,\label{eq:phi-size}\\
 |\Psi(x,z)-\Psi(x,w)|&\lesssim_p L_\Psi|z-w|\,(|z|+|w|)^{p-1},
 \label{eq:spatial-vector-lipschitz}\\
 |\Psi(x,z+w)-\Psi(x,z)-\Psi(x,w)|
 &\lesssim_p L_\Psi\Delta_p(|z|,|w|).\label{eq:pointwise-increment}
\end{align}
Here \eqref{eq:spatial-vector-lipschitz} follows from the polar representation
and the elementary bound
\[
 \min(|z|,|w|)\,\biggl|\frac{z}{|z|}-\frac{w}{|w|}\biggr|\le 2|z-w|
 \qquad(z,w\ne0).
\]
The zero case is immediate. For \eqref{eq:pointwise-increment}, use
the symmetry to assume $|w|\le|z|$. The left-hand side is bounded via
\eqref{eq:phi-size}--\eqref{eq:spatial-vector-lipschitz} by
$C_p L_\Psi|w|\,|z|^{p-1}$. Integrating \eqref{eq:spatial-vector-lipschitz}
and applying H\"older's inequality yields the $L^p$-to-$L^1$ continuity bound
\begin{equation}\label{eq:spatial-vector-Lp-continuity}
 \|\Psi(\cdot,U)-\Psi(\cdot,W)\|_1
 \lesssim_p L_\Psi\|U-W\|_p\bigl(\|U\|_p+\|W\|_p\bigr)^{p-1}.
\end{equation}

Likewise, for any finite vector measure $\nu$ with $\mu=|\nu|$ and total mass
$M=\mu(\R^d)$, both fields $H_n*\nu$ and $H_n(\cdot-a)\nu(\R^d)$ are bounded
in norm by $C_0 M V_n^{-1/p}$. Applying \eqref{eq:spatial-vector-lipschitz}
and integrating the linear translation difference via
\eqref{eq:common-shell-hypotheses} (where $V_n^{-(p-1)/p}$ and $V_n^{1/p'}$
cancel identically) produces the single-shell point-mass replacement bound
\begin{equation}\label{eq:dini-point-mass}
 \bigl\|\Psi(\cdot,H_n*\nu)-\Psi\bigl(\cdot,H_n(\cdot-a)\nu(\R^d)\bigr)\bigr\|_1
 \lesssim L_\Psi M^{p-1}\int_{\R^d}\vartheta\bigl(D_n(y-a)\bigr)\,\dd\mu(y).
\end{equation}

\subsection{Directional slicing, centered transport, and finite shell blocks}

Set $\mu:=|\nu|$ and $\omega_0(t):=\min(1,t)$. For each $n\in\mathbb Z$ and
$x,y\in\R^d$, define the doubled body and local shell masses
\[
 P_n(y):=y+2A_n,\qquad a_n(y):=\mu\bigl(P_n(y)\bigr),\qquad
 b_n(x):=\mu(x-E_n).
\]
The bodies $P_n(y)$ are nested in $n$, and whenever $x\in y+E_n$ we have
$x-A_n\subset y+A_n-A_n=P_n(y)$ by origin-symmetry and convexity of $A_n$;
consequently $b_n(x)\le a_n(y)$ on $y+E_n$.

To control both the mass deficit $a_n(y)-b_n(x)$ and the translation
continuity of small-scale envelopes, we slice convex bodies along the
translation direction. For any convex body $C\subset\R^d$ and any
$h\in\R^d\setminus\{0\}$, Fubini's theorem along lines parallel to $h$ yields
\[
 \|\one_C(\cdot-h)-\one_C\|_1
 =2\int_{h^\perp}\min\bigl(|h|,\ell_C(\xi)\bigr)\,\dd\xi,
\]
where $\ell_C(\xi)$ is the length of the one-dimensional section
$C\cap(\xi+\R h)$. Thus $\|\one_C(\cdot-h)-\one_C\|_1$ is monotone under
inclusion of convex bodies. When $\delta:=D_n(h)<1$, origin-symmetry and
convexity give $(1-\delta)A_n\subset A_n\cap(A_n+h)$, whence
$\|\one_{A_n}(\cdot-h)-\one_{A_n}\|_1\le 2\bigl(1-(1-\delta)^d\bigr)V_n\le 2dV_n\delta$.
For $j\ge n$, combining inclusion monotonicity for $A_j,A_{j+1}\subset A_n$
with the trivial volume bound $2V_j$ establishes the collar estimate
\begin{equation}\label{eq:shell-collar}
 \|\one_{E_j}(\cdot-h)-\one_{E_j}\|_1
 \lesssim \min\bigl\{V_j,\,V_n D_n(h)\bigr\}\qquad(j\ge n,\ h\in\R^d).
\end{equation}

Taking $j=n$ in \eqref{eq:shell-collar} and integrating $a_n(y)-b_n(x)$ over
$x\in y+E_n$ via Tonelli's theorem---noting that each $z\in P_n(y)$ contributes
the collar $(y+E_n)\setminus(z+E_n)$---gives the shell deficit estimate
\begin{equation}\label{eq:shell-deficit}
 \int_{y+E_n}\bigl(a_n(y)-b_n(x)\bigr)\,\dd x
 \lesssim V_n\int_{P_n(y)}\omega_0\bigl(D_n(z-y)\bigr)\,\dd\mu(z).
\end{equation}
For any nonnegative or integrable function $F$ vanishing wherever $b_n=0$,
interchanging the order of integration yields the source-centered
normalization identity
\begin{equation}\label{eq:shell-normalization}
 \int_{\R^d}F(x)\,\dd x
 =\int_{\R^d}\int_{y+E_n}\frac{F(x)}{b_n(x)}\,\dd x\,\dd\mu(y),
\end{equation}
with the integrand set to zero on $\{b_n=0\}$.

The next lemma provides the common summation engine for both the signed path
transfer (Lemma~\ref{lem:path-block-transfer}) and the cross-scale convexity
defect (Proposition~\ref{thm:aniso-interaction}).

\begin{lemma}\label{lem:anchored-transport}
For every finite positive Borel measure $\mu$ on $\R^d$ and every bounded,
nondecreasing, continuous Dini modulus $\omega$ with $\omega(0)=0$,
\begin{equation}\label{eq:centered-transport}
 \sum_{n\in\mathbb Z}\int_{\R^d}a_n(y)^{p-2}
 \int_{P_n(y)}\omega\bigl(D_n(z-y)\bigr)\,\dd\mu(z)\,\dd\mu(y)
 \lesssim \mathcal D_\omega\,\mu(\R^d)^p,
\end{equation}
with zero-mass terms interpreted as zero.
\end{lemma}

\begin{proof}
Fix $y\in\R^d$ and abbreviate $a_n:=a_n(y)$ and $M:=\mu(\R^d)>0$. By
\eqref{eq:common-scale-hypotheses}, the layers $P_k(y)\setminus P_{k+1}(y)$
($k\ge n$) partition $P_n(y)\setminus\{y\}$. Their measures are $a_k-a_{k+1}$,
and on the $k$th layer we have $D_n(z-y)\le 2C_s\,2^{-b_*(k-n)}$; the point
$z=y$ contributes zero because $\omega(0)=0$.

Set $g(t):=\max(M^{p-2},t^{p-2})$ for $t\in(0,M]$. Then $g$ is nonincreasing
on $(0,M]$, and for every $k\ge n$ with $a_k>0$ we have $a_n^{p-2}\le g(a_k)$
(using $a_n\le M$ when $p\ge 2$ and $a_n\ge a_k$ when $1<p<2$). Interchanging
the sums over $n\le k$ and applying \eqref{eq:dini-geometric-sum} gives
\[
 \sum_{n\in\mathbb Z}a_n^{p-2}\int_{P_n(y)}\omega\bigl(D_n(z-y)\bigr)\,\dd\mu(z)
 \lesssim \mathcal D_\omega\sum_{k:a_k>0}g(a_k)(a_k-a_{k+1}).
\]
Since $(a_k)$ is nonincreasing, the intervals $(a_{k+1},a_k]$ are pairwise
disjoint subsets of $(0,M]$, so monotonicity of $g$ yields
\[
 \sum_{k:a_k>0}g(a_k)(a_k-a_{k+1})
 \le \int_0^M g(t)\,\dd t
 \lesssim_p M^{p-1}.
\]
Integrating in $y$ against $\dd\mu(y)$ proves \eqref{eq:centered-transport}.
\end{proof}

We now bound finite decoupled sums of single-shell nonlinear integrals by
$\mathfrak V_\Psi^H$.

\begin{lemma}\label{lem:path-block-transfer}
For every finite $\R^m$-valued Borel measure $\nu$ and all integers $N\le M$,
\begin{equation}\label{eq:path-block-transfer}
 \left|\sum_{n=N}^M\int_{\R^d}\Psi\bigl(x,H_n*\nu(x)\bigr)\,\dd x\right|
 \lesssim \bigl(L_\Psi+\mathfrak V_\Psi^H\bigr)\|\nu\|_{\rm TV}^p.
\end{equation}
\end{lemma}

\begin{proof}
It suffices to consider $\nu\ne0$ and $\mathfrak V_\Psi^H<\infty$.
\emph{Step 1: Point-mass replacement on each shell.}
Fix $n\in[N,M]$ and apply \eqref{eq:shell-normalization} to
$F(x)=\Psi(x,U(x))$ with $U:=H_n*\nu$ (which vanishes wherever $b_n=0$). For
fixed $y\in\R^d$, write $a:=a_n(y)$ and $b:=b_n(x)$, and introduce the
point-mass approximation $U_0(x):=H_n(x-y)\nu(P_n(y))$. On $y+E_n$ we have
$U=H_n*(\nu|_{P_n(y)})$ and $|U|\le C_0 V_n^{-1/p}b$.

Whenever $a>0$, the triangle inequality gives
\[
 \left|\frac{\Psi(x,U)}{b}-\frac{\Psi(x,U_0)}{a}\right|
 \le \frac{|\Psi(x,U)-\Psi(x,U_0)|}{a}+CL_\Psi V_n^{-1}a^{p-2}(a-b).
\]
Indeed, for $b>0$ the denominator difference obeys
\[
 |\Psi(x,U)|(b^{-1}-a^{-1})
 \le C L_\Psi V_n^{-1}b^{p-1}(a-b)/a
 \le C L_\Psi V_n^{-1}a^{p-2}(a-b),
\]
while at $b=0$ the zero-quotient convention reduces the left-hand side to
$|\Psi(x,U_0)|/a$, which is dominated by the first term on the right.
Integrating over $x\in y+E_n$, we control the first term by
\eqref{eq:dini-point-mass} for $\nu|_{P_n(y)}$ (divided by $a$) and the second
by \eqref{eq:shell-deficit}. Summing over $n\in[N,M]$ and integrating in $y$
against $\dd\mu(y)$, Lemma~\ref{lem:anchored-transport} with
$\omega=\vartheta+\omega_0$ bounds the total replacement error by
$CL_\Psi\|\nu\|_{\rm TV}^p$.

\smallskip
\emph{Step 2: Variation bound along nested paths.}
By positive $p$-homogeneity of $\Psi$, the integrated point-mass model equals
$\int_{\R^d}q_{n,y}(v_n(y))\,\dd\mu(y)$, where
$v_n(y):=a_n(y)^{-1/p}\nu(P_n(y))$ (set to zero if $a_n(y)=0$). Since $A_n$ is
closed, $\one_{2A_n}(z-y)$ is jointly Borel measurable, so $y\mapsto a_n(y)$
and $y\mapsto \nu(P_n(y))$ are Borel measurable; likewise, the jointly Borel
integrand of $q_{n,y}(v)$ renders $(y,v)\mapsto q_{n,y}(v)$ Borel measurable,
and every single-shell integral is finite by \eqref{eq:annular-family-size}.

Along the nested path $v_N(y),\ldots,v_M(y)$, if two successive masses satisfy
$a\ge b>0$ with vector masses $u,v\in\R^m$, then $|u-v|\le a-b$ and
$|v|\le b$, so
\[
 |a^{-1/p}u-b^{-1/p}v|
 \le (a-b)a^{-1/p}+b\bigl(b^{-1/p}-a^{-1/p}\bigr)
 \le \Bigl(1+\frac{1}{p}\Bigr)\int_b^a t^{-1/p}\,\dd t.
\]
The bound extends to $b=0$ with zero normalized vector. Telescoping over
$n=N,\ldots,M$ and adding the endpoints $|v_N(y)|,|v_M(y)|\le a_N(y)^{1-1/p}$
gives
\[
 D\bigl(v_N(y),\ldots,v_M(y)\bigr)
 \le C_p\,a_N(y)^{1-1/p}\le C_p\,\|\nu\|_{\rm TV}^{1-1/p}.
\]
Apply \eqref{eq:vector-path-primitive} where the path is nonzero. Zero
paths contribute zero. Integrating against $\dd\mu(y)$ bounds the model by
$C_p\mathfrak V_\Psi^H\|\nu\|_{\rm TV}^p$. Combining Steps~1 and~2 yields
\eqref{eq:path-block-transfer}.
\end{proof}

\section{The convexity defect and cross-scale interaction}
\label{sec:common-interaction}

To pass from decoupled single-shell sums to the full nonlinear integral
$\int\Psi(x,K_H*\nu(x))\,\dd x$, we dominate the cross-scale interaction by a
universal nonnegative convexity defect that relies solely on the shell
geometry $(A_n)$, without any translation modulus on $H_n$.

\begin{proposition}\label{thm:aniso-interaction}
For a finite $\R^m$-valued Borel measure $\nu$, set $\mu:=|\nu|$,
$G_n:=V_n^{-1/p}\one_{E_n}$, and $M_n:=G_n*\mu$. For every finite index set
$I\subset\mathbb Z$, define the nonnegative convexity defect
\[
 \mathcal E_I:=\biggl(\sum_{n\in I}M_n\biggr)^p-\sum_{n\in I}M_n^p.
\]
These defects increase monotonically with $I$, and their pointwise supremum
$\mathcal E:=\sup_{I\text{ finite}}\mathcal E_I$ satisfies
\begin{equation}\label{eq:aniso-interaction}
 \|\mathcal E\|_1\lesssim \|\nu\|_{\rm TV}^p.
\end{equation}
In particular, every finite cross-scale nonlinear remainder obeys
\begin{equation}\label{eq:defect-domination}
 \left|\Psi\biggl(x,\sum_{n\in I}H_n*\nu(x)\biggr)
 -\sum_{n\in I}\Psi\bigl(x,H_n*\nu(x)\bigr)\right|
 \lesssim L_\Psi\,\mathcal E_I(x).
\end{equation}
\end{proposition}

\begin{proof}
\emph{Step 1: $L^q$ translation continuity of the small-scale envelope.}
Set $q:=\max(1,p-1)$, pick $0<\varepsilon<1/q-1/p$, and define the power Dini
modulus $\omega(t):=\min(1,t^\varepsilon)$. Introduce the small-scale envelope
kernel $J_n:=\sum_{j>n}G_j$ and field $W_n:=J_n*\mu$. By geometric volume
decay, the series defining $J_n$ converges in $L^q(\R^d)$. Applying
\eqref{eq:shell-collar} to each $G_j$ ($j>n$) yields
\[
 \|J_n(\cdot-h)-J_n\|_q
 \lesssim \sum_{j>n}V_j^{1/q-1/p}
 \min\biggl(1,\frac{V_n D_n(h)}{V_j}\biggr)^{1/q}
 \lesssim V_n^{1/q-1/p}\omega\bigl(D_n(h)\bigr),
\]
using $\min(1,t)^{1/q}\le t^\varepsilon$ and volume decay (as
$1/q-1/p-\varepsilon>0$) for $D_n(h)<1$, and bounding the minimum by one for
$D_n(h)\ge 1$.

\smallskip
\emph{Step 2: Reverse-anchored shell estimate.}
Fix $y\in\R^d$ and work on $y+E_n$ with normalized measure
$\dd\ell:=V_n^{-1}\dd x$. Set $a:=a_n(y)$, $b:=b_n(x)$, $Z:=V_n^{1/p}W_n$, and
$T:=\int_{P_n(y)}\omega(D_n(z-y))\,\dd\mu(z)\le a$. The kernel $J_n$ vanishes
outside $\operatorname{int}A_{n+1}\subset A_n$, so every source contributing
to $W_n$ on $y+E_n$ lies in $P_n(y)$ and $J_n(x-y)=0$ there. By Young's
inequality $W_n\in L^q$, so almost everywhere on $y+E_n$,
$W_n(x)=\int_{P_n(y)}(J_n(x-z)-J_n(x-y))\,\dd\mu(z)$. Minkowski's integral
inequality together with \eqref{eq:shell-deficit} (and $\omega_0\le\omega$)
gives
\begin{equation}\label{eq:fine-envelopes-on-shell}
 \|Z\|_{L^q(\ell)}+\int_{y+E_n}(a-b)\,\dd\ell\lesssim T.
\end{equation}

For $0<b\le a$ and $z\ge 0$, we bound $\Delta_p(b,z)/b$ as follows:
\begin{itemize}
\item If $1<p<2$ and $b\ge a/2$, then $\Delta_p(b,z)/b\lesssim b^{p-2}z\lesssim a^{p-2}z$.
If $0<b<a/2$, then $\Delta_p(b,z)/b\lesssim z^{p-1}\le a^{p-2}(z+a)$ (splitting
$z\le a$ and $z>a$) and $a\le 2(a-b)$, giving
$\Delta_p(b,z)/b\lesssim a^{p-2}(z+a-b)$ in both cases.
\item If $p\ge 2$, then $\Delta_p(b,z)/b\lesssim a^{p-2}z+z^{p-1}$.
\end{itemize}
Integrating against $\dd\ell$ via \eqref{eq:fine-envelopes-on-shell}, using
$\ell(y+E_n)\le 1$ and $T^{p-1}\le a^{p-2}T$ for $p\ge 2$, and noting that
$M_n=V_n^{-1/p}b$ and $W_n=V_n^{-1/p}Z$ on $y+E_n$, we obtain
\begin{equation}\label{eq:reverse-source-anchor}
 \int_{y+E_n}\frac{\Delta_p(M_n,W_n)}{b}\,\dd x
 =\int_{y+E_n}\frac{\Delta_p(b,Z)}{b}\,\dd\ell
 \lesssim a^{p-2}T.
\end{equation}
The quotients vanish when $b=0$, and the estimate is trivial if $a=0$.
Applying \eqref{eq:shell-normalization} and summing via
Lemma~\ref{lem:anchored-transport} gives
$\sum_{n\in\mathbb Z}\int\Delta_p(M_n,W_n)\,\dd x\lesssim\mu(\R^d)^p$.

\smallskip
\emph{Step 3: Telescoping and defect domination.}
On any symmetric block $I=[-J,J]$, telescoping in decreasing index order
yields
\[
 \mathcal E_{[-J,J]}
 =\sum_{|n|\le J}\Delta_p\biggl(M_n,\sum_{n<j\le J}M_j\biggr).
\]
Extending the summands by zero outside $[-J,J]$ and letting $J\to\infty$ by
monotone convergence gives $\mathcal E=\sum_{n\in\mathbb Z}\Delta_p(M_n,W_n)$,
proving \eqref{eq:aniso-interaction}. Because $|H_n*\nu|\le C_0 M_n$,
telescoping over any finite set $I$ in decreasing index order and applying
\eqref{eq:pointwise-increment} establishes \eqref{eq:defect-domination}.
\end{proof}

\section{Path realization and proof of the general classification}
\label{sec:proof-sharp}

We now assemble the tools from Sections~\ref{sec:pointmass-energy}
and~\ref{sec:common-interaction} to prove
Theorem~\ref{thm:path-family-classification} and
Corollary~\ref{thm:annular-family-classification}, and then record their
immediate measure-theoretic corollaries for concentrated sources and absolute
estimates.

\subsection{The decoupled series identity and large-scale tails}

For a finite $\R^m$-valued Borel measure $\nu$, put $U_n:=H_n*\nu$. Since
$\sum_{n\le 0}\|U_n\|_\infty+\sum_{n>0}\|U_n\|_1<\infty$, both $\sum_n U_n$
and---by \eqref{eq:phi-size}---$\sum_n\Psi(\cdot,U_n)$ converge absolutely almost
everywhere. Defining $\mathcal E_I$ by monotone exhaustion for arbitrary
$I\subset\mathbb Z$ and passing to the limit in \eqref{eq:defect-domination}
yields
\begin{equation}\label{eq:spatial-vector-mixed-remainder}
 \left|\Psi\biggl(x,\sum_{n\in I}U_n(x)\biggr)
 -\sum_{n\in I}\Psi\bigl(x,U_n(x)\bigr)\right|
 \lesssim L_\Psi\,\mathcal E_I(x).
\end{equation}
Taking $I=\mathbb Z$ gives the decoupled representation
\begin{equation}\label{eq:spatial-vector-density-series}
 \Psi(\cdot,K_H*\nu)=\sum_{n\in\mathbb Z}\Psi(\cdot,U_n)+R_\nu,
 \qquad \|R_\nu\|_1\lesssim L_\Psi\|\nu\|_{\rm TV}^p.
\end{equation}
The same argument applies to $|\Psi|$, since $L_{|\Psi|}\le L_\Psi$.
In particular, for every finite vector measure $\nu$,
\begin{equation}\label{eq:absolute-density-decoupling}
 \biggl\|\,|\Psi(\cdot,K_H*\nu)|-\sum_{n\in\mathbb Z}|\Psi(\cdot,U_n)|\,\biggr\|_1
 \lesssim L_\Psi\|\nu\|_{\rm TV}^p.
\end{equation}
Thus absolute integrability is equivalent to
$\sum_n\|\Psi(\cdot,U_n)\|_1<\infty$, by Tonelli's theorem.

\begin{lemma}
\label{lem:spatial-vector-endpoints}
Let $\nu$ be a compactly supported finite $\R^m$-valued Borel measure with
total mass $v=\nu(\R^d)$.
\begin{enumerate}
\item Membership in $\mathcal X_H$ is characterized by
\begin{equation}\label{eq:compact-source-fine-Lp}
 \nu\in\mathcal X_H\quad\Longleftrightarrow\quad
 \sum_{n>0}\|U_n\|_p^p<\infty.
\end{equation}
In particular, $L_c^\infty(\R^d;\R^m)\subset\mathcal X_H$.
\item The integrand $\Psi(\cdot,K_H*\nu)$ belongs to $L^1(\R^d)$ if and only
if $\mathfrak T_\Psi^H(a,v)<\infty$ and $\sum_{n>0}\|\Psi(\cdot,U_n)\|_1<\infty$.
Finiteness of $\mathfrak T_\Psi^H(a,v)$ is independent of the center $a\in\R^d$
and of any fixed scale cutoff.
\item If $\nu\in\mathcal X_H$ and $\mathfrak T_\Psi^H(a,v)<\infty$, then
$\sum_n\Psi(\cdot,U_n)$ converges absolutely in $L^1(\R^d)$, and
$\Psi(\cdot,H_{\le N}*\nu)$ converges in $L^1(\R^d)$ to $0$ as $N\to-\infty$
and to $\Psi(\cdot,K_H*\nu)$ as $N\to+\infty$. Moreover, mass-preserving
smooth mollification $\nu_\delta=\nu*\chi_\delta$ yields
$\Psi(\cdot,K_H*\nu_\delta)\to\Psi(\cdot,K_H*\nu)$ in $L^1(\R^d)$. In
particular, every mean-zero $\nu\in\mathcal X_H$ satisfies
$K_H*\nu\in L^p(\R^d;\R^\ell)$.
\end{enumerate}
\end{lemma}

\begin{proof}
\emph{Proof of \textup{(i)}.}
Since $H_{\le 0}*\nu$ is bounded and $\sum_{n>0}U_n$ has compact support,
$K_H*\nu\in L^p_{\mathrm{loc}}$ is equivalent to $\sum_{n>0}U_n\in L^p(\R^d)$.
Applying \eqref{eq:spatial-vector-mixed-remainder} with $\Psi(x,z)=|z|^p$ on
$I=\{n>0\}$ proves \eqref{eq:compact-source-fine-Lp}. For
$f\in L_c^\infty(\R^d;\R^m)$, Young's inequality and
$\sum_{n>0}\|H_n\|_1^p<\infty$ verify this condition.

\smallskip
\emph{Proof of \textup{(ii)}.}
The coarse shellwise differences between $\Psi(\cdot,U_n)$ and
$\Psi(\cdot,H_n(\cdot-a)v)$ are summable in $L^1$ by \eqref{eq:dini-point-mass}
and the Dini sum, since $\nu$ is compactly supported. The criterion now
follows from \eqref{eq:absolute-density-decoupling}. Applying the same
shellwise comparison to $\nu=v\delta_b$ proves center independence, while
changing the cutoff only adds or removes finitely many integrable shells.

\smallskip
\emph{Large-scale point-mass comparison.}
Fix a compact set $K$ containing $\supp\nu$ and $a$. For $n\le 0$, geometric
contraction gives $D_n(z-y)\le C_K\,2^{b_*n}=:\delta_n$ for all $z,y\in K$.
Let $\omega(t)=\min(1,t^\varepsilon)$ be the power modulus from
Proposition~\ref{thm:aniso-interaction}. Then for every $N\le 0$,
\begin{equation}\label{eq:uniform-coarse-atomic-comparison}
 \bigl\|\Psi(\cdot,H_{\le N}*\nu)-\Psi\bigl(\cdot,H_{\le N}(\cdot-a)v\bigr)\bigr\|_1
 \lesssim L_\Psi\|\nu\|_{\rm TV}^p
 \sum_{n\le N}\bigl[\omega(\delta_n)+\vartheta(\delta_n)\bigr].
\end{equation}
Indeed, for $n\le N$ and $y\in\supp\nu\subset K$, the transport quantity in
\eqref{eq:reverse-source-anchor} obeys $T\le a_n(y)\omega(\delta_n)$, so
\eqref{eq:shell-normalization} bounds the $n$th integrated reverse defect by
$C\|\nu\|_{\rm TV}^p\omega(\delta_n)$. Summing over $n\le N$ controls
$\|\mathcal E_{(-\infty,N]}\|_1$, and replacing each $\Psi(\cdot,U_n)$ by
$\Psi(\cdot,H_n(\cdot-a)v)$ via \eqref{eq:dini-point-mass} (using
disjointness of the shells $a+E_n$) proves
\eqref{eq:uniform-coarse-atomic-comparison}; the right-hand side converges by
\eqref{eq:dini-geometric-sum}.

\smallskip
\emph{Proof of \textup{(iii)}.}
Assume $\nu\in\mathcal X_H$ and $\mathfrak T_\Psi^H(a,v)<\infty$. By
\textup{(i)}, the fine terms are absolutely summable in $L^1$, so the
full density is integrable by \textup{(ii)} and all shellwise terms are
absolutely summable in $L^1$ by \eqref{eq:absolute-density-decoupling}.
As $N\to+\infty$, the shellwise sum
converges in $L^1$ by absolute summability and the remainder by domination by
$\mathcal E\in L^1$; as $N\to-\infty$, $L^1$ convergence to $0$ follows from
\eqref{eq:uniform-coarse-atomic-comparison}. Applying \textup{(ii)} to
$\Psi(x,z)=|z|^p$ when $v=0$ gives $K_H*\nu\in L^p(\R^d;\R^\ell)$.

Finally, let $\chi_\delta$ ($0<\delta\le 1$) be a nonnegative smooth radial
probability mollifier supported in $B(0,\delta)$, and set
$\nu_\delta:=\nu*\chi_\delta$. Then $\|\nu_\delta\|_{\rm TV}\le\|\nu\|_{\rm TV}$,
every $\nu_\delta$ has mass $v$, and all $\supp\nu_\delta$ lie in a common
compact set $K$. By \eqref{eq:mollification-associativity} and
$K_H*\nu\in L^p_{\mathrm{loc}}$, we have $K_H*\nu_\delta\to K_H*\nu$ in
$L^p_{\mathrm{loc}}$ as $\delta\downarrow 0$; for any fixed $N\le 0$,
\eqref{eq:spatial-vector-Lp-continuity} yields $L^1$ convergence of the
nonlinear integrands on the compact set $K+A_{N+1}$. Outside $K+A_{N+1}$, only
$H_{\le N}*\nu_\delta$ and $H_{\le N}*\nu$ remain, and their $L^1$ nonlinear
norms vanish uniformly in $\delta\in(0,1]$ as $N\to-\infty$ by
\eqref{eq:uniform-coarse-atomic-comparison} and $\mathfrak T_\Psi^H(a,v)<\infty$.
\end{proof}

\subsection{Realization of finite mass paths}

\begin{lemma}\label{lem:path-realization}
Given a center $a\in\R^d$ and a finite path $v_N,\ldots,v_M\in\R^m$, extend it
by zero outside $[N,M]$ and set $c_t:=v_t-v_{t+1}$,
$L:=\sum_{t\in\mathbb Z}|c_t|=D(v_N,\ldots,v_M)$, and
\[
 f:=\sum_{t\in\mathbb Z}c_t\eta_t^a,\qquad
 \eta_t^a:=V_t^{-1}\one_{a+A_t}.
\]
Then $f\in L_c^\infty(\R^d;\R^m)$ satisfies $\int_{\R^d}f=0$ and
$\|f\|_1\le L$. Moreover,
\begin{align}
 \sum_{n\in\mathbb Z}\bigl\|\Psi(\cdot,H_n*f)-\Psi\bigl(\cdot,H_n(\cdot-a)v_n\bigr)\bigr\|_1
 &\lesssim L_\Psi L^p,\label{eq:path-realization-error}\\
 \left|\int_{\R^d}\Psi\bigl(x,K_H*f(x)\bigr)\,\dd x-\sum_{n=N}^M q_{n,a}(v_n)\right|
 &\lesssim L_\Psi L^p.\label{eq:path-realization-full}
\end{align}
\end{lemma}

\begin{proof}
Telescoping gives $v_n=\sum_{t\ge n}c_t$ and $\sum_t c_t=0$. Each $\eta_t^a$
has unit integral and $\|\eta_t^a\|_p=V_t^{-1/p'}$. Split $f=f_{\ge n}+f_{<n}$
according to $t\ge n$ and $t<n$. For $t\ge n$, every $x\in a+A_t$ satisfies
$D_n(x-a)\le C_s\,2^{-b_*(t-n)}$, so \eqref{eq:dini-point-mass} compares
$\Psi(\cdot,H_n*f_{\ge n})$ with $\Psi(\cdot,H_n(\cdot-a)v_n)$. For $t<n$,
Young's inequality and \eqref{eq:annular-family-size} give
$\|(H_n c_t)*\eta_t^a\|_p\lesssim |c_t|(V_n/V_t)^{1/p'}$. Since both $H_n*f$
and $H_n*f_{\ge n}$ have $L^p$ norms at most $CL$,
\eqref{eq:spatial-vector-Lp-continuity} bounds their nonlinear difference.
Combining both estimates yields
\[
 \bigl\|\Psi(\cdot,H_n*f)-\Psi\bigl(\cdot,H_n(\cdot-a)v_n\bigr)\bigr\|_1
 \lesssim L_\Psi L^{p-1}\biggl[
 \sum_{t\ge n}|c_t|\,\vartheta\bigl(C_s\,2^{-b_*(t-n)}\bigr)
 +\sum_{t<n}|c_t|\,(V_n/V_t)^{1/p'}\biggr].
\]
Summing in $n\in\mathbb Z$ via \eqref{eq:dini-geometric-sum} and geometric
volume decay proves \eqref{eq:path-realization-error}, and adding
\eqref{eq:spatial-vector-density-series} gives \eqref{eq:path-realization-full}.
\end{proof}

\subsection{Proof of the master classification}

\begin{proof}[Proof of Theorem~\ref{thm:path-family-classification}]
By Lemma~\ref{lem:spatial-vector-endpoints}, the integrand
$\Psi(\cdot,K_H*\nu)$ for $\nu\in\mathcal X_H$ belongs to $L^1(\R^d)$ if and
only if its total mass $v=\nu(\R^d)$ satisfies $\mathfrak T_\Psi^H(a_0,v)<\infty$.
Mean-zero sources always qualify, and every $v\in\R^m$ arises as the mass of a
smooth compactly supported source. In either integrable branch, the
mollification convergence in Lemma~\ref{lem:spatial-vector-endpoints}\textup{(iii)}
preserves every admissible constant on $C_c^\infty(\R^d;\R^m)$, so the optimal
constants on $\mathcal X_H$ and $C_c^\infty(\R^d;\R^m)$ coincide.

For sufficiency, let $\nu\in\mathcal X_H$ lie in the appropriate integrable
branch. Applying Lemma~\ref{lem:path-block-transfer} to finite blocks $[N,M]$,
passing to the limit $N\to-\infty$, $M\to+\infty$ via absolute $L^1$
summability of $\sum_n\Psi(\cdot,U_n)$, and adding
\eqref{eq:spatial-vector-density-series} yields
the upper bound $C_{\Psi,0}^{H,*}\lesssim L_\Psi+\mathfrak V_\Psi^H$ in the
mean-zero branch. Under \eqref{eq:vector-all-direction-tails}, the same
argument gives $C_\Psi^{H,*}\lesssim L_\Psi+\mathfrak V_\Psi^H$ for all masses.
Conversely, testing the mean-zero sources constructed in
Lemma~\ref{lem:path-realization} gives
$\mathfrak V_\Psi^H\le C_{\Psi,0}^{H,*}+CL_\Psi$, and $C_{\Psi,0}^{H,*}\le
C_\Psi^{H,*}$ in the full-source branch. These inequalities
establish \eqref{eq:vector-zero-mean-quantitative} and
\eqref{eq:vector-full-quantitative}.
\end{proof}

\subsection{Spherical path reduction and proof of the scalar classification}

\begin{lemma}\label{lem:scalar-path-reduction}
For $0<\delta\le 1$, let $\mathfrak P_\delta$ denote the supremum of
$|\sum_{n=N}^M q_{n,a}(u_n)|$ over all centers $a\in\R^d$, integer blocks
$N\le M$, and unit-vector paths $(u_n)_{n=N}^M\subset\mathbb S^{m-1}$ with
interior variation $\sum_{n=N}^{M-1}|u_{n+1}-u_n|\le\delta$. Then
\begin{equation}\label{eq:short-spherical-testing}
 \mathfrak P_\delta\le(2+\delta)^p\mathfrak V_\Psi^H,
 \qquad
 \mathfrak V_\Psi^H\le\bigl(2^{-p}+2^{1-p}/\delta\bigr)\mathfrak P_\delta.
\end{equation}
In the scalar case $m=1$, $\mathfrak P_\delta$ equals
$\mathfrak P:=\sup_{a,\sigma=\pm 1,N\le M}|\sum_{n=N}^M q_{n,a}(\sigma)|$, and
\begin{equation}\label{eq:scalar-sharp-testing}
 \mathfrak V_\Psi^H=2^{-p}\mathfrak P.
\end{equation}
\end{lemma}

\begin{proof}
Any unit-vector path with interior variation at most $\delta$ has
zero-extended variation $D\le 2+\delta$, yielding the first inequality in
\eqref{eq:short-spherical-testing}. For the converse, take any nonzero finite
path $v_n=r_nu_n$ padded by zeros, with $r_n=|v_n|$, $|u_n|=1$ (chosen
arbitrarily where $r_n=0$), $R:=\max_n r_n$, and $D\ge2R$.
By positive $p$-homogeneity and the layer-cake
representation,
\[
 \sum_{n=N}^M q_{n,a}(v_n)
 =\int_0^R ph^{p-1}\sum_{n:\,r_n>h}q_{n,a}(u_n)\,\dd h.
\]
At each height $h\in(0,R)$, greedily partition every maximal consecutive
component of $\{n:r_n>h\}$ into maximal successive unit paths of variation at
most $\delta$. Each cut consumes strictly more than $\delta$ of angular
variation (counting the preceding segment together with the discarded jump).
If $N(h)$ counts the active components and $L(h)$ denotes their total interior
angular variation, then at most $N(h)+L(h)/\delta$ segments arise, so
\[
 \left|\sum_{n=N}^M q_{n,a}(v_n)\right|
 \le \mathfrak P_\delta\int_0^R ph^{p-1}\bigl(N(h)+L(h)/\delta\bigr)\,\dd h.
\]
Since $N(h)\ge 1$ on $(0,R)$ and
$\int_0^R N(h)\,\dd h=D(r_N,\ldots,r_M)/2\le D/2$, convexity of $t\mapsto t^p$
yields
\[
 \int_0^R ph^{p-1}N(h)\,\dd h
 \le R^p+pR^{p-1}(D/2-R)\le (D/2)^p=2^{-p}D^p.
\]
For the angular term, the identity $|ru-sw|^2=(r-s)^2+rs|u-w|^2$ ($r,s>0$,
$|u|=|w|=1$) gives $\min(r,s)|u-w|\le |ru-sw|$, whence
\[
 \int_0^R ph^{p-1}L(h)\,\dd h
 =\sum_n\min(r_n,r_{n+1})^p|u_{n+1}-u_n|
 \le R^{p-1}D\le 2^{1-p}D^p,
\]
proving \eqref{eq:short-spherical-testing}.

When $m=1$, any path on $\mathbb S^0=\{-1,1\}$ of variation at most
$\delta\le 1$ is constant, so $\mathfrak P_\delta=\mathfrak P$. Moreover, at
height $h\in(0,R)$ the active segments are the consecutive components of
$\{v_n>h\}$ and $\{v_n<-h\}$; their combined count $N_\pm(h)$ satisfies
$N_\pm(h)\ge 1$ and $\int_0^R N_\pm(h)\,\dd h=D/2$, giving
$\mathfrak V_\Psi^H\le 2^{-p}\mathfrak P$. Taking the supremum over constant
signed paths proves the reverse inequality.
\end{proof}

\smallskip
\noindent\emph{Sharpness of the variation dependence.}
The order $\delta^{-1}$ in \eqref{eq:short-spherical-testing} is necessary.
Identify $\R^2$ with $\mathbb C$, take $d=m=\ell=2$, and use
$A_n=\overline{B(0,2^{-n})}$ and
$H_n(x)v=(-1)^nV_n^{-1/p}\one_{E_n}(x)v$. These kernels obey uniform size
and linear translation bounds by \eqref{eq:shell-collar}. For odd $k\ge3$
let $\Phi_k(z)=k^{-1}|z|^p\cos(k\arg z)$ for $z\ne0$ and $\Phi_k(0)=0$.
The spherical norms $L_{\Phi_k}$ are uniformly bounded, since angular
distance on $\mathbb S^1$ is at most $\pi/2$ times chordal distance.
Oddness of $\Phi_k$ and $|E_n|/V_n=3/4$ give the unit-vector responses
\[
 q_{n,a}(e^{i\theta})=\frac{3(-1)^n}{4k}\cos(k\theta).
\]
The alternating signs have uniformly bounded partial sums. Abel summation
therefore bounds every consecutive unit-path response by $C(k^{-1}+\delta)$
for interior variation at most $\delta$. A one-shell
path gives $\mathfrak P_{1/k}\asymp k^{-1}$. In contrast, the $k$-term path
starting at $n=0$ and alternating between $1$ and $e^{i\pi/k}$ has
$D\le2+\pi$ and response $3/4$. Thus $\mathfrak V_{\Phi_k}^H\gtrsim1$,
and the upper bound in \eqref{eq:short-spherical-testing} gives
$\mathfrak V_{\Phi_k}^H\asymp1$. All structural data are independent of $k$.

\begin{proof}[Proof of Corollary~\ref{thm:annular-family-classification}]
Apply Theorem~\ref{thm:path-family-classification} and
Lemma~\ref{lem:scalar-path-reduction} to $\Psi=b\Phi$, for which
$L_\Psi=BL_\Phi$ and $\mathfrak V_\Psi^H=2^{-p}\mathfrak P_b^H$. Scalar
homogeneity reduces \eqref{eq:vector-all-direction-tails} to the two signs in
\eqref{eq:family-absolute-tail}. In either applicable branch, the optimal
constant $C^*$ obeys
\[
 C^*\lesssim BL_\Phi+\mathfrak P_b^H,
 \qquad \mathfrak P_b^H\le 2^p C^*+CBL_\Phi,
\]
This proves the comparisons with $BL_\Phi$. Multiplying the bounds by
at most $\max(1,L_\Phi)$ gives \eqref{eq:family-zero-mean-quantitative}
and \eqref{eq:family-finite-quantitative}. When $b\in L^1\cap L^\infty$, the
coarse tail is bounded by $CL_\Phi\|b\|_1\sum_{n\le 0}V_n^{-1}<\infty$. For
$b\equiv 1$, translation invariance identifies $\mathfrak P_1^H$ with
$\Lambda_H$.
\end{proof}

\subsection{Concentrated measure sources and absolute estimates}
\label{sec:concentrated-measures}

We conclude the general convex-shell theory with two immediate measure-level
consequences of Theorem~\ref{thm:path-family-classification} and
Corollary~\ref{thm:annular-family-classification}. Set
$h_n(x):=\|H_n(x)\|_{\rm op}$ and define the small-scale energy of a finite
positive Borel measure $\mu$ by
\begin{equation}\label{eq:fine-source-energy}
 \mathscr E_H(\mu):=\sum_{n\ge 0}\|h_n*\mu\|_p^p.
\end{equation}

\begin{lemma}\label{lem:fine-measure-envelope}
We have $\mathscr E_H(\mu)<\infty$ if and only if
$\sum_{n\ge 0}h_n*\mu\in L^p(\R^d)$, in which case the small-scale tails
converge to $0$ in $L^p(\R^d)$. This holds for every compactly supported $L^q$
density ($q>1$) and whenever
\begin{equation}\label{eq:source-mass-modulus-condition}
 \sum_{n\ge 0}\omega_n^{p-1}<\infty,
 \qquad \omega_n:=\sup_{a\in\R^d}\mu(a+A_n).
\end{equation}
For diagonal shells, \eqref{eq:source-mass-modulus-condition} includes every
Frostman measure satisfying $\sup_a\mu(a+\{\rho_{\mathbf a}<r\})\le Ar^s$
for some $s\in(0,d]$.
\end{lemma}

\begin{proof}
Applying \eqref{eq:spatial-vector-mixed-remainder} to $h_n$ and
$\Psi(x,z)=|z|^p$ yields the equivalence, and dominated convergence gives the
tail limit. For $g\in L_c^q$ ($q>1$), putting $q_0:=\min(q,p)$ and applying
Young's inequality gives $\|h_n*|g|\|_p\lesssim V_n^{1-1/q_0}\|g\|_{q_0}$,
whose $p$th powers decay geometrically. In general,
$\|h_n*\mu\|_\infty\lesssim V_n^{-1/p}\omega_n$ and
$\|h_n*\mu\|_1\lesssim V_n^{1/p'}\mu(\R^d)$, so interpolation gives
$\|h_n*\mu\|_p^p\lesssim\mu(\R^d)\omega_n^{p-1}$. For diagonal shells,
letting $r\downarrow 2^{-n}$ gives $\omega_n\le A\,2^{-sn}$.
\end{proof}

\begin{corollary}
\label{cor:concentrated-measure-sources}
Every compactly supported finite $\R^m$-valued Borel measure $\nu$ with
$\mathscr E_H(|\nu|)<\infty$ belongs to $\mathcal X_H$. Both branches of
Theorem~\ref{thm:path-family-classification} hold on this class with the exact
same optimal constants as on $C_c^\infty(\R^d;\R^m)$.
\end{corollary}

\begin{proof}
Because $|H_n*\nu|\le h_n*|\nu|$,
Lemma~\ref{lem:spatial-vector-endpoints}\textup{(i)} gives $\nu\in\mathcal X_H$.
Since $C_c^\infty(\R^d;\R^m)$ lies in this class, the optimal constants
coincide by Theorem~\ref{thm:path-family-classification}.
\end{proof}

Among compactly supported finite measures, the energy condition can be
strictly stronger than membership in $\mathcal X_H$. On isotropic shells,
take $H_n=V_n^{-1/p}\one_{E_n}A$ with $A\ne0$ and $Av=0$ for some $v\ne0$.
Then $v\delta_0\in\mathcal X_H$ because its potential vanishes, whereas
$\mathscr E_H(|v|\delta_0)=\infty$.

For unsigned (absolute) estimates with scalar sources, no local $L^p$
regularity is needed at all:

\begin{corollary}
\label{cor:absolute-finite-measure-sources}
Assume the hypotheses of Corollary~\ref{thm:annular-family-classification} and
set
\begin{equation}\label{eq:absolute-atomic-measure-test}
 \mathfrak A_b^{H,{\rm abs}}
 :=\sup_{a\in\R^d,\ \sigma=\pm 1}\sum_{n\in\mathbb Z}\int_{\R^d}
 |b(x)|\,\bigl|\Phi\bigl(\sigma H_n(x-a)\bigr)\bigr|\,\dd x.
\end{equation}
Then the absolute estimate
\begin{equation}\label{eq:absolute-finite-measure-estimate}
 \int_{\R^d}|b(x)|\,\bigl|\Phi(K_H*\nu(x))\bigr|\,\dd x
 \le C_{\rm abs}\|\nu\|_{\rm TV}^p
\end{equation}
holds for all finite real signed Borel measures $\nu$ if and only if
$\mathfrak A_b^{H,{\rm abs}}<\infty$. Its optimal constant $C_{b,{\rm abs}}^{H,*}$
satisfies $B+C_{b,{\rm abs}}^{H,*}\asymp B+\mathfrak A_b^{H,{\rm abs}}$ and is
unchanged upon restricting to $L_c^\infty(\R^d;\R)$.
\end{corollary}

\begin{proof}
Applying Corollary~\ref{thm:annular-family-classification} to $|b|$ and
$|\Phi|$ gives $\mathfrak P_{|b|}^H=\mathfrak A_b^{H,{\rm abs}}$, which also
dominates the coarse tail, establishing \eqref{eq:absolute-finite-measure-estimate}
on $L_c^\infty$. For an arbitrary finite signed measure $\nu$, put
$f_j:=(\nu|_{B(0,j)})*\chi_{1/j}$ for a smooth probability mollifier supported
in $B(0,1/j)$. Because $K_H\in L^\infty+L^1$ and
$\|\nu|_{B(0,j)}-\nu\|_{\rm TV}\to 0$, we have
$K_H*(\nu|_{B(0,j)})\to K_H*\nu$ in $L^1_{\mathrm{loc}}$, and mollification
preserves local $L^1$ convergence via \eqref{eq:mollification-associativity}.
Passing to an almost-everywhere convergent subsequence and applying Fatou's
lemma preserves every admissible $L_c^\infty$ constant. Testing
$\nu=\sigma\delta_a$ yields $\mathfrak A_b^{H,{\rm abs}}\le C_{b,{\rm abs}}^{H,*}$
by shell disjointness.
\end{proof}

\section{Anisotropic dilations and measurable targets}
\label{sec:anisotropic-measurable}

We now specialize the general classification of Section~\ref{sec:proof-sharp}
to homogeneous kernels under positive diagonal dilations, general expanding
matrix dilations, and Euclidean measurable targets, proving
Theorems~\ref{thm:anisotropic-mazya-phi}
and~\ref{thm:sharp-measurable-characterization}, recovering
$C^{1,\mathrm{Dini}}$ boundaries, and establishing universal target rigidity.

\subsection{Anisotropic polar coordinates and sharp-annular regularity}
\label{sec:polar-geometry}

For a fixed diagonal weight vector $\mathbf a\in(0,\infty)^d$ with
$\sum_{j=1}^d a_j=d$, the critical exponents satisfy $p\beta=d$ and
$(p-1)\beta=\alpha$. Any unnormalized positive diagonal dilation with weights
$(b_j)$ and homogeneous dimension $Q=\sum_{j=1}^d b_j$ is normalized by
setting $a_j=db_j/Q$ and $\alpha=d\alpha_{\mathbf b}/Q$, preserving $p$ and
scaling the polar flux weight by $d/Q$.

Let $\mathcal K:\R^d\setminus\{0\}\to\R^\ell$ satisfy
$\mathcal K(\delta_t^{\mathbf a}x)=t^{-\beta}\mathcal K(x)$ for $t>0$, with
$\mathcal K(0)=0$ and bounded Borel angular restriction
$\widetilde{\mathcal K}=\mathcal K|_{\mathbb S^{d-1}}$. Define the flux
weight and sharp dyadic annular shells
\begin{equation}\label{eq:aniso-weight-shells}
 w_{\mathbf a}(\omega):=\sum_{j=1}^d a_j\omega_j^2,
 \qquad
 \mathcal K_n(x):=\mathcal K(x)\,\one_{\{2^{-n-1}\le\rho_{\mathbf a}(x)<2^{-n}\}}.
\end{equation}

\begin{lemma}\label{lem:aniso-polar}
For every nonnegative measurable function $G$ on $\R^d$,
\begin{equation}\label{eq:aniso-polar}
 \int_{\R^d}G(x)\,\dd x
 =\int_0^\infty\int_{\mathbb S^{d-1}}
 G(\delta_r^{\mathbf a}\omega)\,w_{\mathbf a}(\omega)\,\dd\sigma(\omega)\,r^{d-1}\dd r.
\end{equation}
Consequently, for $\sigma\in\{-1,1\}$ and every $n\in\mathbb Z$,
\begin{equation}\label{eq:aniso-shell-mass}
 \int_{\R^d}\Phi(\sigma\mathcal K_n)\,\dd x=(\log 2)\,A_\sigma^{\mathbf a},
 \qquad
 A_\sigma^{\mathbf a}:=\int_{\mathbb S^{d-1}}w_{\mathbf a}(\omega)\,
 \Phi\bigl(\sigma\widetilde{\mathcal K}(\omega)\bigr)\,\dd\sigma(\omega).
\end{equation}
\end{lemma}

\begin{proof}
The map $(r,\omega)\mapsto\delta_r^{\mathbf a}\omega$ is a global $C^1$
diffeomorphism from $(0,\infty)\times\mathbb S^{d-1}$ onto $\R^d\setminus\{0\}$.
Factoring $\delta_r^{\mathbf a}$ from its differential yields
$\det(\delta_r^{\mathbf a})=r^d$; the remaining radial column is
$r^{-1}(a_j\omega_j)_{j=1}^d$, whose inner product with the outward unit
normal $\omega$ to $\mathbb S^{d-1}$ is $r^{-1}w_{\mathbf a}(\omega)$, proving
\eqref{eq:aniso-polar}. In the nonlinear shell integral, the radial density
becomes $r^{d-1-p\beta}\dd r=\dd r/r$, and absolute integrability on one shell
gives \eqref{eq:aniso-shell-mass}.
\end{proof}

By homogeneity, $\mathcal K_n(x)=2^{\beta n}\mathcal K_0(\delta_{2^n}^{\mathbf a}x)$,
so any unit-shell $L^1$ translation bound
\begin{equation}\label{eq:rough-kernel-dini}
 \|\mathcal K_0(\cdot-h)-\mathcal K_0\|_1
 \le C_\vartheta\,\vartheta(|h|)\qquad(h\in\R^d)
\end{equation}
with $\vartheta$ as in \eqref{eq:dini-modulus} scales directly to
\begin{equation}\label{eq:aniso-dini-translation}
 \|\mathcal K_n(\cdot-h)-\mathcal K_n\|_1
 \le C_\vartheta\,2^{-\alpha n}\vartheta\bigl(|\delta_{2^n}^{\mathbf a}h|\bigr),
\end{equation}
rendering $H_n=\mathcal K_n$ admissible in \eqref{eq:annular-family-hypotheses}
(and analogously in operator norm for matrix-valued kernels). The following
lemma records broad sufficient angular regularity criteria for
\eqref{eq:rough-kernel-dini}.

\begin{lemma}
\label{cor:angular-dini-sufficient}
The unit-shell translation bound \eqref{eq:rough-kernel-dini} holds whenever
any of the following conditions is satisfied:
\begin{enumerate}
\item the bounded angular restriction $\widetilde{\mathcal K}$ has a uniform
Dini modulus $\omega$ satisfying \eqref{eq:dini-modulus} (in particular,
$\widetilde{\mathcal K}\in C^\gamma(\mathbb S^{d-1})$ for some
$0<\gamma\le 1$);
\item there is a finite union $\Gamma\subset\mathbb S^{d-1}$ of compact $C^1$
embedded $(d-2)$-dimensional submanifolds (with or without boundary) such that
the restrictions of $\widetilde{\mathcal K}$ to the connected components of
$\mathbb S^{d-1}\setminus\Gamma$ share a common Lipschitz constant;
\item the bounded unit shell $\mathcal K_0$ belongs to $BV(\R^d;\R^\ell)$.
\end{enumerate}
In cases \textup{(ii)} and \textup{(iii)}, and for Lipschitz angular parts in
\textup{(i)}, one may take $\vartheta(t)=\min(t,1)$.
\end{lemma}

\begin{proof}
Both $\rho_{\mathbf a}$ and the angular projection
$\pi_{\mathbf a}(x)=\delta_{\rho_{\mathbf a}(x)^{-1}}^{\mathbf a}x$ are smooth
on a neighborhood of the unit annulus $\{1/2\le\rho_{\mathbf a}<1\}$. In case
\textup{(i)}, on the overlap with a small translate the pointwise difference
is $O(|h|+\omega(C|h|))$, while the symmetric difference of the two annuli has
volume $O(|h|)$. Combining this with the trivial $L^1$ bound for large $h$
gives \eqref{eq:rough-kernel-dini} with the continuous Dini modulus
$\vartheta(t)=\min(t,1)+\omega(Ct)$.

In case \textup{(ii)}, covering $\Gamma$ by finitely many Lipschitz graph
charts encloses its spherical $\epsilon$-neighborhood in strips of thickness
$O(\epsilon)$, yielding spherical measure $O(\epsilon)$. Along a small
translation segment $[x-h,x]$ inside a compact neighborhood of the unit
annulus, $\pi_{\mathbf a}$ moves by at most $C|h|$. Whenever $\pi_{\mathbf a}(x)$
lies outside the $C|h|$-neighborhood of $\Gamma$, the segment image remains
within a single connected component of $\mathbb S^{d-1}\setminus\Gamma$, where
the difference is $O(|h|)$. By Lemma~\ref{lem:aniso-polar}, the volume of the
exceptional directions is $O(|h|)$, giving a linear modulus.

In case \textup{(iii)}, smoothing and the fundamental theorem of calculus give
$\|\mathcal K_0(\cdot-h)-\mathcal K_0\|_1\le |h|\,|D\mathcal K_0|(\R^d)$.
\end{proof}

\subsection{Homogeneous kernels, weighted targets, and anisotropic cancellation}
\label{sec:anisotropic-mazya-phi}

Assuming \eqref{eq:rough-kernel-dini} and putting $H_n=\mathcal K_n$, we write
$\mathfrak P_{b,\mathbf a}:=\mathfrak P_b^H$ for $b\in L^\infty(\R^d;\R)$ with
$B=\|b\|_\infty$, and define the continuous annular testing seminorm and
large-scale absolute tail by
\begin{align}
 \mathfrak P_{b,\mathbf a}^{\rm c}
 &:=\sup_{\substack{c\in\R^d,\ \sigma=\pm 1\\0<r<R<\infty}}
 \left|\int_{\{r<\rho_{\mathbf a}(x-c)<R\}}
 b(x)\Phi\bigl(\sigma\mathcal K(x-c)\bigr)\,\dd x\right|,
 \label{eq:weighted-continuous-cocycle}\\
 \mathfrak T_{b,\mathbf a}(c_0)
 &:=\max_{\sigma=\pm 1}\int_{\{\rho_{\mathbf a}(x-c_0)\ge 1\}}
 |b(x)|\,\bigl|\Phi\bigl(\sigma\mathcal K(x-c_0)\bigr)\bigr|\,\dd x.
 \label{eq:aniso-target-absolute-tail}
\end{align}
When $b=\one_E$, we replace the subscript $b$ by $E$.

\Needspace*{3\baselineskip}
\begin{lemma}
\label{lem:continuous-dyadic-cocycle}
For every $b\in L^\infty(\R^d;\R)$,
\begin{equation}\label{eq:continuous-dyadic-equivalence}
 \mathfrak P_{b,\mathbf a}
 \le \mathfrak P_{b,\mathbf a}^{\rm c}
 \le \mathfrak P_{b,\mathbf a}+CB,
\end{equation}
where one can take
\[
 C=2\log 2\max_{\sigma=\pm1}\int_{\mathbb S^{d-1}}w_{\mathbf a}(\omega)
        |\Phi(\sigma\widetilde{\mathcal K}(\omega))|\,\dd\sigma(\omega).
\]
The boundary convention on shells is immaterial. For a fixed geometric ratio
$\varrho>1$ the same comparison holds with $\log\varrho$ in place of $\log 2$.
\end{lemma}

\begin{proof}
Every dyadic block is a continuous annulus up to null boundaries, giving the
first inequality. Conversely, removing at most two partial dyadic shells from
$\{r<\rho_{\mathbf a}(x-c)<R\}$ leaves a (possibly empty) dyadic block. By
Lemma~\ref{lem:aniso-polar}, each partial shell costs at most $CB/2$.
\end{proof}

\begin{corollary}
\label{cor:weighted-homogeneous-target}
Assume \eqref{eq:rough-kernel-dini} and let $b\in L^\infty(\R^d;\R)$ with
$B=\|b\|_\infty$.
\begin{enumerate}
\item For every mean-zero $f\in L_c^\infty(\R^d;\R)$, the integral in
\begin{equation}\label{eq:aniso-zero-mean-target-estimate}
 \left|\int_{\R^d}b(x)\Phi(\mathcal K*f(x))\,\dd x\right|
 \le C_{b,0}\|f\|_1^p
\end{equation}
converges absolutely. The estimate holds uniformly on this class if and only
if $\mathfrak P_{b,\mathbf a}^{\rm c}<\infty$, with optimal constant
$C_{b,0,\mathbf a}^*$ satisfying
$B+C_{b,0,\mathbf a}^*\asymp B+\mathfrak P_{b,\mathbf a}^{\rm c}$.
\item Absolute convergence for every $f\in L_c^\infty(\R^d;\R)$ is equivalent
to $\mathfrak T_{b,\mathbf a}(c_0)<\infty$. Under this condition, the
unconstrained estimate
\begin{equation}\label{eq:weighted-homogeneous-estimate}
 \left|\int_{\R^d}b(x)\Phi(\mathcal K*f(x))\,\dd x\right|
 \le C_b\|f\|_1^p
 \qquad\bigl(f\in L_c^\infty(\R^d;\R)\bigr)
\end{equation}
holds if and only if $\mathfrak P_{b,\mathbf a}^{\rm c}<\infty$, with optimal
constant $C_{b,\mathbf a}^*$ satisfying
\begin{equation}\label{eq:weighted-homogeneous-comparison}
 B+C_{b,\mathbf a}^*\asymp B+\mathfrak P_{b,\mathbf a}^{\rm c}.
\end{equation}
\end{enumerate}
Finiteness of $\mathfrak T_{b,\mathbf a}(c_0)$ is independent of $c_0$ and of
the radial cutoff, and is automatic whenever $b\in L^1\cap L^\infty$ (in
particular for $b=\one_E$ with $|E|<\infty$). The comparison constants depend
only on $d$, $\ell$, $\mathbf a$, $\alpha$, $\|\widetilde{\mathcal K}\|_\infty$,
the data in \eqref{eq:rough-kernel-dini}, and $L_\Phi$.
\end{corollary}

\begin{proof}
By \eqref{eq:aniso-dini-translation}, the shells $H_n=\mathcal K_n$ satisfy
\eqref{eq:annular-family-hypotheses}. Applying
Corollary~\ref{thm:annular-family-classification} together with
\eqref{eq:continuous-dyadic-equivalence} yields all claims, as the discrete
coarse tail $\mathfrak T_b^H(c_0)$ differs from
\eqref{eq:aniso-target-absolute-tail} by at most one bounded integrable shell.
\end{proof}

\begin{proof}[Proof of Theorem~\ref{thm:anisotropic-mazya-phi}]
Because $\widetilde{\mathcal K}$ is Lipschitz,
Lemma~\ref{cor:angular-dini-sufficient} and \eqref{eq:aniso-dini-translation}
show that $H_n=\mathcal K_n$ satisfies \eqref{eq:annular-family-hypotheses}.
By Lemma~\ref{lem:aniso-polar}, each whole-space shell total equals
$(\log 2)A_\sigma^{\mathbf a}$ independently of $n$, so $\Lambda_H<\infty$ if
and only if $A_\sigma^{\mathbf a}=0$ for both signs $\sigma=\pm 1$.
Corollary~\ref{thm:annular-family-classification} with $b\equiv 1$ completes
the proof.
\end{proof}

\begin{remark}
\label{rem:weighted-full-sphere}
If $b\in L^\infty(\R^d;\R)$ satisfies $|\{b\ne0\}|>0$ and
$\mathfrak P_{b,\mathbf a}^{\rm c}<\infty$, then both weighted full-sphere
coefficients $A_\sigma^{\mathbf a}$ in \eqref{eq:aniso-shell-mass} must
vanish, without any integrability assumption on $b$. Indeed, for fixed
$\sigma=\pm 1$, set $G_n(z):=\Phi(\sigma\mathcal K_n(z))$,
$\kappa:=\int G_n=(\log 2)A_\sigma^{\mathbf a}$, and
$Q_n(c):=\int G_n(z)b(c+z)\,\dd z$. The kernels $G_n$ have constant $L^1$ norm
and supports shrinking to the origin, so on every ball $B_R$,
\[
 \|Q_n-\kappa b\|_{L^1(B_R)}
 \le \|G_n\|_1\sup_{z\in\supp G_n}\|b(\cdot+z)-b\|_{L^1(B_R)}\longrightarrow 0
 \qquad(n\to+\infty)
\]
by $L^1_{\mathrm{loc}}$ continuity of translations. Consequently the Ces\`aro
averages $M^{-1}\sum_{n=1}^M Q_n$ converge to $\kappa b$ in $L^1(B_R)$, while
$|\sum_{n=1}^M Q_n|\le\mathfrak P_{b,\mathbf a}^{\rm c}$ forces them to
converge uniformly to $0$. Since $b\not\equiv 0$, we conclude $\kappa=0$.
\end{remark}

\subsection{Expanding matrix dilations in Lyapunov sections}
\label{sec:semisimple-dilations}

The convex-shell theory extends directly from diagonal dilations to arbitrary
expanding real matrix groups $T_t=t^A=\exp((\log t)A)$.

\begin{corollary}\label{cor:semisimple-dilations}
Let $A\in\R^{d\times d}$ have $\operatorname{Re}\lambda>0$ for every
eigenvalue $\lambda$ and $\operatorname{tr}A=d$, and set $T_t:=t^A$ ($t>0$).
Choose a symmetric positive definite matrix $Q$ such that $A^TQ+QA>0$ (for
example, $Q=\int_0^\infty e^{-sA^T}e^{-sA}\,\dd s$), write
$|x|_Q:=(x^TQx)^{1/2}$, $B_Q:=\{|x|_Q<1\}$, and $\Sigma:=\partial B_Q$.
Define $\rho_A(x)$ for $x\ne 0$ as the unique positive solution of
$|T_{1/\rho_A(x)}x|_Q=1$, with $\rho_A(0)=0$.

Let $\beta=d/p$, $\alpha=d-\beta$, and let $\Phi:\R^\ell\to\R$ be positively
$p$-homogeneous with $\Phi(0)=0$ and $L_\Phi<\infty$. Suppose
$\mathcal K:\R^d\setminus\{0\}\to\R^\ell$ satisfies
$\mathcal K(T_tx)=t^{-\beta}\mathcal K(x)$ ($t>0$), $\mathcal K(0)=0$, and
$\mathcal K|_\Sigma$ is bounded with a Dini modulus of continuity. Then for
$H_n:=\mathcal K\one_{\{2^{-n-1}\le\rho_A<2^{-n}\}}$, all conclusions of
Corollary~\ref{thm:annular-family-classification} hold for every
$b\in L^\infty(\R^d;\R)$, with the continuous annular test and tail
\eqref{eq:weighted-continuous-cocycle}--\eqref{eq:aniso-target-absolute-tail}
formed using $\rho_A$. In particular, the whole-space mean-zero estimate holds
if and only if
\begin{equation}\label{eq:semisimple-cancellation}
 \int_\Sigma (Az\cdot n_\Sigma(z))\,\Phi\bigl(\sigma\mathcal K(z)\bigr)\,
 \dd\mathcal H^{d-1}(z)=0
 \qquad(\sigma=\pm 1),
\end{equation}
where $n_\Sigma(z)=Qz/|Qz|$ is the exterior unit normal to $\Sigma$.
\end{corollary}

\begin{proof}
Because every eigenvalue of $A$ has positive real part, $e^{-sA}$ decays
exponentially as $s\to\infty$, so $Q=\int_0^\infty e^{-sA^T}e^{-sA}\,\dd s$ is
symmetric positive definite and solves $A^TQ+QA=I$. For any symmetric positive
definite $Q$ with $A^TQ+QA>0$, differentiating
\[
 \frac{\dd}{\dd s}|e^{sA}x|_Q^2=(e^{sA}x)^T(A^TQ+QA)e^{sA}x
 \ge 2c\,|e^{sA}x|_Q^2
\]
yields $|T_tx|_Q\le t^c|x|_Q$ for $0<t\le 1$ and $|T_tx|_Q\ge t^c|x|_Q$ for
$t\ge 1$. Hence $t\mapsto|T_tx|_Q$ increases strictly from $0$ to $\infty$ for
each $x\ne 0$; by the implicit function theorem, $\rho_A$ is smooth on
$\R^d\setminus\{0\}$ and satisfies $\rho_A(T_tx)=t\rho_A(x)$. The ellipsoids
$A_n:=T_{2^{-n}}\overline{B_Q}=\{\rho_A\le 2^{-n}\}$ are nested
origin-symmetric convex bodies obeying $A_{n+k}\subset 2^{-ck}A_n$, with
volumes $V_n=2^{-dn}|B_Q|$ and Minkowski norms $D_n(h)=|T_{2^n}h|_Q$.

On the unit shell $\{1/2\le\rho_A<1\}$, both $\rho_A$ and the projection
$x\mapsto T_{1/\rho_A(x)}x\in\Sigma$ are smooth. Arguing as in
Lemma~\ref{cor:angular-dini-sufficient} and using \eqref{eq:shell-collar}
gives a bounded Dini modulus $\vartheta$ for $H_0$ in the $Q$-norm. Exact
scaling $H_n(x)=2^{\beta n}H_0(T_{2^n}x)$ and $\det(T_{2^n})=2^{dn}$ then
verify \eqref{eq:common-shell-hypotheses}, so
Corollary~\ref{thm:annular-family-classification} applies.

Finally, $(r,z)\mapsto T_rz$ is a global diffeomorphism from
$(0,\infty)\times\Sigma$ onto $\R^d\setminus\{0\}$. Factoring $T_r$ from its
differential and projecting $r^{-1}Az$ onto $n_\Sigma(z)=Qz/|Qz|$ yields
\[
 \dd x=r^{d-1}(Az\cdot n_\Sigma(z))\,\dd r\,\dd\mathcal H^{d-1}(z),
 \qquad Az\cdot n_\Sigma(z)=\frac{z^T(A^TQ+QA)z}{2|Qz|}>0.
\]
Critical homogeneity yields the radial measure $\dd r/r$, so the whole-space
shell totals equal $\log 2$ times the integrals in
\eqref{eq:semisimple-cancellation}, and partial-shell errors remain $O(B)$ as
in Lemma~\ref{lem:continuous-dyadic-cocycle}.
\end{proof}

When $A$ is semisimple in the form $S^{-1}AS=D+\Omega$ with
$D=\operatorname{diag}(a_j)>0$, $\Omega^T=-\Omega$, and $[D,\Omega]=0$,
choosing $Q=S^{-T}S^{-1}$ gives $\Sigma=S\mathbb S^{d-1}$,
$\rho_A(Sy)=\rho_{\mathbf a}(y)$, and flux element
$|\det S|\,w_{\mathbf a}(\omega)\,\dd\sigma(\omega)$. If the symmetric part of
$A$ is positive definite, one may simply take $Q=I$.

\subsection{Euclidean logarithmic profiles and Dini-regular boundaries}
\label{sec:cocycle-necessity}

Returning to the Euclidean setting of
Theorem~\ref{thm:sharp-measurable-characterization} ($\mathbf a=\mathbf 1=(1,\ldots,1)$),
set
\begin{equation}\label{eq:annular-kernels}
 K_n(x):=K(x)\,\one_{\{2^{-n-1}\le|x|<2^{-n}\}},
 \qquad
 b_n^E(a,\epsilon):=\int_E\Phi\bigl(\epsilon K_n(x-a)\bigr)\,\dd x.
\end{equation}
Write $\mathfrak P_E:=\mathfrak P_{\one_E}^H$ and
$\mathfrak P_E^{\rm cont}:=\mathfrak P_{\one_E,\mathbf 1}^{\rm c}$ for
$H_n=K_n$. When replacing $(K,\Phi)$ by another Euclidean pair
$(\mathcal K,\Psi)$, we append $[\mathcal K,\Psi]$ in brackets after the
functional or after the arguments of $b_n^E$.

\begin{proof}[Proof of Theorem~\ref{thm:sharp-measurable-characterization}]
Apply Corollary~\ref{cor:weighted-homogeneous-target}\textup{(ii)} with
$\mathbf a=\mathbf 1$ and $b=\one_E$. Because $0<|E|<\infty$, we have $B=1$
and the large-scale tail $\mathfrak T_{E,\mathbf 1}(c_0)$ is automatically
finite. The quantitative equivalence \eqref{eq:main-continuous-quantitative}
follows from \eqref{eq:weighted-homogeneous-comparison} and
\eqref{eq:continuous-dyadic-equivalence}.
\end{proof}

To analyze the geometric content of $\mathfrak P_E^{\rm cont}$, put
$\psi_\epsilon(\omega):=\Phi(\epsilon\widetilde K(\omega))$ and define the
angular occupancy profile
\begin{equation}\label{eq:angular-occupancy}
 \mathcal A_E(a,r,\epsilon):=\int_{\mathbb S^{d-1}}
 \one_{E_{\mathrm B}}(a+r\omega)\psi_\epsilon(\omega)\,\dd\sigma(\omega),
\end{equation}
where $E_{\mathrm B}$ is any Borel representative of $E$; by polar Fubini,
$\mathcal A_E(a,r,\epsilon)$ is independent of $E_{\mathrm B}$ for
$\dd r/r$-almost every $r>0$.

\begin{proposition}
\label{prop:logarithmic-shell-formula}
For every measurable $E\subset\R^d$, center $a\in\R^d$, sign $\epsilon=\pm 1$,
and integer $n\in\mathbb Z$,
\begin{equation}\label{eq:logarithmic-shell-formula}
 b_n^E(a,\epsilon)
 =\int_{2^{-n-1}}^{2^{-n}}\mathcal A_E(a,r,\epsilon)\,\frac{\dd r}{r}.
\end{equation}
\end{proposition}

\begin{proof}
Pass to polar coordinates centered at $a$, observing that
$r^{d-1-p(d-\alpha)}\dd r=\dd r/r$.
\end{proof}

\begin{remark}
\label{prop:bounded-logarithmic-primitive}
Setting $q_{a,\epsilon}(t):=\mathcal A_E(a,e^{-t},\epsilon)$, we have
\[
 \mathfrak P_E^{\rm cont}
 =\sup_{a\in\R^d,\ \epsilon=\pm 1,\ s<t}
 \left|\int_s^t q_{a,\epsilon}(u)\,\dd u\right|.
\]
For each bounded measurable profile $q$ on $\R$, this interval supremum equals
$2\inf_{B'=q}\|B\|_\infty$, the infimum ranging over bounded distributional
primitives of $q$ (namely, the oscillation of $\int_0^t q(u)\,\dd u$ centered
by an additive constant).
\end{remark}

We now show that Theorem~\ref{thm:sharp-measurable-characterization} recovers---and
extends to $C^{1,\mathrm{Dini}}$ boundaries---Stolyarov's tangent-hemisphere
criterion for smooth domains~\cite[Theorem~1.2]{StolyarovDomains}. For a
boundary point $\xi$ of a $C^1$ open set $\Omega$ with exterior unit normal
$\nu(\xi)$, let $H_\xi^-:=\{\omega\in\mathbb S^{d-1}:\omega\cdot\nu(\xi)<0\}$.

\begin{corollary}\label{thm:smooth-domain-recovery}
Let $\Omega\subset\R^d$ be a nonempty bounded open set with $C^1$ boundary,
lying locally on one side of each boundary graph. Suppose there exist $r_0>0$
and a continuous nondecreasing modulus $\omega:[0,r_0]\to[0,\infty)$ with
$\omega(0)=0$ and $\int_0^{r_0}\omega(t)\,\dd t/t<\infty$ such that
\begin{equation}\label{eq:c1beta-tangent-error}
 \bigl|(\Omega\mathbin\triangle\mathcal H_\xi)\cap B(\xi,r)\bigr|
 \le r^d\omega(r),
 \qquad \mathcal H_\xi:=\{x:(x-\xi)\cdot\nu(\xi)<0\},
\end{equation}
uniformly for $\xi\in\partial\Omega$ and $0<r<r_0$. Then
$\mathfrak P_\Omega<\infty$ if and only if
\begin{equation}\label{eq:smooth-hemisphere}
 \int_{H_\xi^-}\psi_\epsilon\,\dd\sigma=0
 \qquad(\xi\in\partial\Omega,\ \epsilon=\pm 1).
\end{equation}
These identities are equivalent to cancellation over every hemisphere of
$\mathbb S^{d-1}$ and imply full-sphere cancellation
$\int_{\mathbb S^{d-1}}\psi_\epsilon\,\dd\sigma=0$. Their necessity holds for
all bounded one-sided $C^1$ open sets. When the identities and
\eqref{eq:c1beta-tangent-error} with a Dini modulus hold,
\[
 \sup_{a,\epsilon}\int_0^\infty
 \bigl|\mathcal A_\Omega(a,r,\epsilon)\bigr|\,\frac{\dd r}{r}
 \le C_{d,K,\Phi}\biggl(1+\int_0^{r_0}\omega(t)\,\frac{\dd t}{t}
 +r_0^{-d}|\Omega|\biggr),
\]
and the same bound controls $\sup_{a,\epsilon}\sum_n|b_n^\Omega(a,\epsilon)|$.
In particular, this criterion covers all bounded $C^{1,\mathrm{Dini}}$ open
sets and recovers the $C^{1,\beta_\partial}$ theorem
\cite[Theorem~1.2]{StolyarovDomains} for $0<\beta_\partial\le 1$.
\end{corollary}

\begin{proof}
For every unit vector $v\in\mathbb S^{d-1}$, a maximum of $x\mapsto x\cdot v$
on $\overline\Omega$ occurs at a boundary point with exterior normal $v$.
Thus \eqref{eq:smooth-hemisphere} yields cancellation over every hemisphere,
and summing complementary hemispheres yields full-sphere cancellation.

\smallskip
\emph{Necessity.}
At each $\xi\in\partial\Omega$, differentiability of the boundary graph gives
$|(\Omega\mathbin\triangle\mathcal H_\xi)\cap B(\xi,r)|=o(r^d)$ as $r\downarrow0$.
Using this estimate in the shell of outer radius $2^{-n}$ gives
$b_n^\Omega(\xi,\epsilon)=(\log 2)\int_{H_\xi^-}\psi_\epsilon\,\dd\sigma+o(1)$
as $n\to\infty$. A nonzero hemisphere integral would force consecutive block
sums to diverge.

\smallskip
\emph{Sufficiency.}
Fix $a\in\R^d$, put $\tau:=\dist(a,\partial\Omega)$, and write $r_n:=2^{-n}$.
We split the dyadic scales into four ranges:
\begin{enumerate}
\item $r_n>r_0/4$: these shells lie in $\{|x-a|>r_0/8\}$, so polar Fubini
bounds their total contribution by $Cr_0^{-d}|\Omega|$;
\item $r_n\le\tau$ (with $r_n\le r_0/4$): the ball $B(a,r_n)$ lies entirely
inside or outside $\Omega$, so $\mathcal A_\Omega(a,r,\epsilon)=0$ by
full-sphere cancellation;
\item $\tau<r_n<2\tau$ (with $r_n\le r_0/4$): at most one dyadic scale lies
in this range, contributing $O(1)$;
\item $2\tau\le r_n\le r_0/4$: choose a nearest point $\xi\in\partial\Omega$
and let $\mathcal H=\{x:(x-a)\cdot\nu(\xi)<0\}$ be the translate of
$\mathcal H_\xi$ centered at $a$. Since $\mathcal A_{\mathcal H}(a,r,\epsilon)=0$,
$B(a,r_n)\subset B(\xi,2r_n)$, and $\mathcal H_\xi\mathbin\triangle\mathcal H$
is a slab of width at most $\tau$, polar Fubini and
\eqref{eq:c1beta-tangent-error} give
\[
 \int_{r_n/2}^{r_n}\bigl|\mathcal A_\Omega(a,r,\epsilon)\bigr|\,\frac{\dd r}{r}
 \lesssim r_n^{-d}\bigl|(\Omega\mathbin\triangle\mathcal H)\cap B(a,r_n)\bigr|
 \lesssim \omega(2r_n)+\frac{\tau}{r_n}.
\]
Bound $\omega(2r_n)$ by its logarithmic average on $[2r_n,4r_n]$. These
intervals have disjoint interiors in $[0,r_0]$, so their sum is bounded by
the Dini integral times $1/\log 2$. The geometric sum of $\tau/r_n$ is at most
one for $\tau>0$, and all its terms vanish when $\tau=0$.
\end{enumerate}
The dyadic sum bound follows from \eqref{eq:logarithmic-shell-formula}. For a
$C^{1,\mathrm{Dini}}$ boundary, recentering local graphs gives
$|h(y')-h(0)-Dh(0)y'|\le|y'|\omega_1(|y'|)$, and integrating the resulting
slab verifies \eqref{eq:c1beta-tangent-error} with $\omega(r)=C\omega_1(Cr)$
on a smaller uniform chart radius.
\end{proof}

\subsection{Universal target rigidity}

We conclude Section~\ref{sec:anisotropic-measurable} by characterizing when a
kernel--nonlinearity pair satisfies the critical inequality on \emph{every}
measurable target.

\begin{corollary}
\label{prop:family-universal-targets}
Under \eqref{eq:annular-family-hypotheses} with $m=1$, define
\[
 S_{\rm fine}:=\max_{\sigma=\pm 1}\sum_{n\ge 0}\|\Phi(\sigma H_n)\|_1,
 \qquad
 S_{\rm all}:=\max_{\sigma=\pm 1}\sum_{n\in\mathbb Z}\|\Phi(\sigma H_n)\|_1,
\]
and abbreviate $\mathfrak P_E^H:=\mathfrak P_{\one_E}^H$,
$C_E^{H,*}:=C_{\one_E}^{H,*}$, and $C_{E,0}^{H,*}:=C_{\one_E,0}^{H,*}$.
\begin{enumerate}
\item Every measurable set $E$ of finite measure admits
\eqref{eq:family-finite-target-estimate} (with $b=\one_E$ and a constant
depending on $E$) if and only if $S_{\rm fine}<\infty$. In that case,
\[
 C_E^{H,*}\le C\bigl[1+S_{\rm fine}+\log(1+|E|)\bigr].
\]
\item Every measurable set $E\subset\R^d$ (including $|E|=\infty$) admits the
mean-zero estimate on $C_c^\infty(\R^d;\R)$ with a constant depending on $E$
if and only if $S_{\rm all}<\infty$. This is also equivalent to a uniform
constant over all finite-measure targets and all $f\in L_c^\infty(\R^d;\R)$,
and
\begin{equation}\label{eq:family-universal-quantitative}
 1+\sup_{E\text{ measurable}}C_{E,0}^{H,*}
 \asymp 1+S_{\rm all}
 \asymp 1+\sup_{|E|<\infty}C_E^{H,*}.
\end{equation}
\item Whenever $S_{\rm all}<\infty$, the absolute whole-space estimate
\begin{equation}\label{eq:family-universal-absolute}
 \int_{\R^d}|\Phi(K_H*f)|\,\dd x
 \le C(1+S_{\rm all})\|f\|_1^p
 \qquad\bigl(f\in L_c^\infty(\R^d;\R)\bigr)
\end{equation}
holds, so all measurable targets share a common constant without the mean-zero
restriction.
\end{enumerate}
\end{corollary}

\begin{proof}
For $n<0$, the absolute value of each shell integral over $E$ is at most
$C\min(1,|E|\,2^{dn})$. Summing and adding the fine-shell contribution gives
$\mathfrak P_E^H\le S_{\rm fine}+C\log(1+|E|)$. The scalar classification
then gives sufficiency in \textup{(i)}.

For necessity, set $J=\{n\ge 0\}$ in \textup{(i)} and $J=\mathbb Z$ in
\textup{(ii)}. Since the shells have disjoint supports, the sign sets
$E_{\sigma,\pm}^J:=\bigcup_{n\in J}\{\pm\Phi(\sigma H_n)>0\}$ satisfy
\[
 \sum_{n\in J}\|\Phi(\sigma H_n)\|_1
 \le \mathfrak P_{E_{\sigma,+}^J}^H+\mathfrak P_{E_{\sigma,-}^J}^H
\]
at center $a=0$, with $|E_{\sigma,\pm}^{\{n\ge 0\}}|<\infty$.
Corollary~\ref{thm:annular-family-classification} gives necessity in both
\textup{(i)} and \textup{(ii)}, together with the lower bound for
$\sup_E C_{E,0}^{H,*}$.

By Corollary~\ref{cor:absolute-finite-measure-sources} (with $b\equiv 1$), the
optimal constant in \eqref{eq:family-universal-absolute} is comparable to
$S_{\rm all}$ after adding $1$. For each $f\in L_c^\infty$, taking the
supremum of $|\int_E\Phi(K_H*f)|$ over finite-measure sets $E$ recovers the
larger of the positive and negative integrals of $\Phi(K_H*f)$, which lies
between half and one times $\int|\Phi(K_H*f)|$. This establishes
\eqref{eq:family-universal-quantitative}.
\end{proof}

\begin{corollary}
\label{thm:universal-target-classification}
Let $\mathcal K$ be homogeneous of degree $-\beta$ under $\delta_t^{\mathbf a}$
with bounded Borel angular restriction $\widetilde{\mathcal K}$, and assume
$H_n=\mathcal K_n$ satisfies \eqref{eq:annular-family-hypotheses}. Let
$\Phi:\R^\ell\to\R$ be positively $p$-homogeneous with $\Phi(0)=0$ and
$L_\Phi<\infty$. The following are equivalent:
\begin{enumerate}
\item for every measurable $E\subset\R^d$ there is a finite $C_E$ such that
every $f\in L_c^\infty(\R^d;\R)$ satisfies $\Phi(\mathcal K*f)\in L^1(E)$ and
$|\int_E\Phi(\mathcal K*f)|\le C_E\|f\|_1^p$;
\item the two signed angular profiles vanish almost everywhere:
\begin{equation}\label{eq:pointwise-angular-annihilation}
 \Phi\bigl(\widetilde{\mathcal K}(\omega)\bigr)
 =\Phi\bigl(-\widetilde{\mathcal K}(\omega)\bigr)=0
 \qquad\text{for almost every }\omega\in\mathbb S^{d-1}.
\end{equation}
\end{enumerate}
The equivalence persists if targets are restricted to finite measure and/or
sources to mean zero. Under \textup{(ii)}, every target has vanishing
point-mass functional and $\int_{\R^d}|\Phi(\mathcal K*f)|\,\dd x\le C\|f\|_1^p$
for all $f\in L_c^\infty(\R^d;\R)$. For continuous $\widetilde{\mathcal K}$,
the almost-everywhere vanishing in \textup{(ii)} is equivalent to pointwise
vanishing.
\end{corollary}

\begin{proof}
Lemma~\ref{lem:aniso-polar} gives the constant shell norm
\[
 \|\Phi(\sigma\mathcal K_n)\|_1
 = (\log 2)\int_{\mathbb S^{d-1}}w_{\mathbf a}|\Phi(\sigma\widetilde{\mathcal K})|\,\dd\sigma.
\]
Since $w_{\mathbf a}\ge\min_j a_j>0$, condition \textup{(ii)} is equivalent
to $S_{\rm fine}<\infty$ and implies $S_{\rm all}=0$. Thus
Corollary~\ref{prop:family-universal-targets} gives the absolute estimate and
\textup{(i)}. Conversely, even on finite-measure targets with mean-zero
sources, \textup{(i)} implies the signed test is finite for each such target.
The tail is automatic there, so the scalar classification gives the
full-source estimate and Corollary~\ref{prop:family-universal-targets} forces
$S_{\rm fine}<\infty$. For continuous angular restrictions, the profiles
are continuous and hence vanish everywhere if they vanish almost everywhere.
\end{proof}

\section{Higher-order elliptic PDEs, vector kernels, and nonstationary shells}
\label{sec:pde-and-extensions}

\subsection{Homogeneous constant-coefficient elliptic operators}

Let $P(\partial)=\sum_{|\gamma|=m}p_\gamma\partial^\gamma$ be a real
constant-coefficient homogeneous differential operator of even order
$m=2q\ge 2$ whose Fourier symbol obeys the ellipticity bound
\begin{equation}\label{eq:elliptic-symbol}
 \pi_P(\xi):=\sum_{|\gamma|=m}p_\gamma(i\xi)^\gamma
 \ge c_P|\xi|^m\qquad(\xi\in\R^d,\ c_P>0);
\end{equation}
a negative elliptic symbol is reduced to \eqref{eq:elliptic-symbol} by
replacing $P$ with $-P$. Let $L(\partial)$ be a real constant-coefficient
homogeneous vector-valued differential operator of order $k\ge 0$ taking
values in $\R^\ell$, and assume
\begin{equation}\label{eq:elliptic-critical-orders}
 \alpha:=m-k\in(0,d),\qquad \beta:=d-\alpha,\qquad
 p:=\frac{d}{\beta},\qquad \nu:=m-d.
\end{equation}

By the standard construction of elliptic fundamental solutions~\cite{Camus},
there exist a real-valued locally integrable function $\Gamma$, smooth on
$\R^d\setminus\{0\}$, and a real homogeneous polynomial $Q$ (of degree $\nu$
if $\nu\ge 0$, and $Q\equiv 0$ if $\nu<0$) satisfying in $\mathcal D'(\R^d)$
\begin{equation}\label{eq:elliptic-associated-homogeneity}
 P\Gamma=\delta_0,\qquad
 \Gamma(t\cdot)=t^\nu\bigl(\Gamma+(\log t)Q\bigr)\qquad(t>0).
\end{equation}
Indeed, set $F(z):=\mathcal F^{-1}(\pi_P^z)$. Because $\pi_P$ is real and
even, $F(\bar z)=\overline{F(z)}$ where $\pi_P^z$ is locally integrable and
hence by meromorphic continuation, so all Laurent coefficients at $z=-1$ are
real. Polar Taylor expansion shows that $F(z)$ has at most a simple pole at
$z=-1$ with Fourier-side residue supported at the origin. Comparing Laurent
coefficients in $F(z)(t\cdot)=t^{-d-mz}F(z)$ yields
\eqref{eq:elliptic-associated-homogeneity} with $\Gamma$ equal to the finite
part at $z=-1$ and $Q=-m\operatorname{Res}_{z=-1}F(z)$; the identity
$PF(z)=\mathcal F^{-1}(\pi_P^{z+1})$ gives $P\Gamma=\delta_0$. Elliptic
regularity yields smoothness off the origin, and the scaling law bounds the
off-origin representative by $C|x|^\nu(1+|\log|x||)$. Its locally integrable
extension obeys the same scaling law. The difference from $\Gamma$ is
supported at zero and homogeneous of degree $\nu>-d$, so it vanishes.

\begin{corollary}
\label{cor:elliptic-measurable-target}
Under \eqref{eq:elliptic-symbol}--\eqref{eq:elliptic-associated-homogeneity},
set $K:=L\Gamma$ on $\R^d\setminus\{0\}$ and $K(0):=0$. Let
$\Phi:\R^\ell\to\R$ be positively $p$-homogeneous with $\Phi(0)=0$ and
$L_\Phi<\infty$. For $b\in L^\infty(\R^d;\R)$ with $B=\|b\|_\infty$, define
$\mathfrak P_b[K,\Phi]:=\mathfrak P_b^H$ and
$\mathfrak P_b^{\rm c}[K,\Phi]:=\mathfrak P_{b,\mathbf 1}^{\rm c}$ using the
isotropic shells $H_n=K_n$. Then $\mathfrak P_b[K,\Phi]<\infty$ if and only if
there exists $C<\infty$ such that
\begin{equation}\label{eq:elliptic-measurable-target}
 \left|\int_{\R^d}b(x)\Phi\bigl(Lu(x)\bigr)\,\dd x\right|
 \le C\|Pu\|_{L^1(\R^d)}^p
 \qquad\bigl(u\in C_c^\infty(\R^d;\R)\bigr).
\end{equation}
Every integral on the left converges absolutely, and the optimal constant
$C_b^{P,L,*}$ coincides with the optimal constant in
\eqref{eq:aniso-zero-mean-target-estimate} for $(K,\Phi)$. In particular,
\begin{equation}\label{eq:elliptic-quantitative-comparison}
 B+C_b^{P,L,*}\asymp_{d,P,L,\Phi} B+\mathfrak P_b[K,\Phi]
 \asymp_{d,P,L,\Phi} B+\mathfrak P_b^{\rm c}[K,\Phi],
\end{equation}
with comparison constants independent of $b$.
\end{corollary}

\begin{proof}
Since $m>m-d=\nu$ and $k=m-\alpha>m-d=\nu$, both $PQ$ and $LQ$ vanish. By
\eqref{eq:elliptic-associated-homogeneity}, $L\Gamma$ is homogeneous of degree
$\alpha-d\in(-d,0)$ with smooth spherical restriction; because no nonzero
distribution supported at the origin has degree $\alpha-d>-d$, $K=L\Gamma$ in
$\mathcal D'(\R^d)$.

For every $u\in C_c^\infty(\R^d;\R)$, $f:=Pu$ is smooth, compactly supported,
and mean zero, and $K*f=Lu$; thus $C_b^{P,L,*}\le C_{b,0,\mathbf 1}^*$.
Conversely, given any mean-zero $f\in C_c^\infty(\R^d;\R)$, set
$v:=\Gamma*f$ and $u_R:=\chi(\cdot/R)v$ for a cutoff $\chi\in C_c^\infty(B_2)$
with $\chi\equiv 1$ on $B_1$. Subtract $\nabla^j\Gamma(x)$ in the convolution
using $\int f=0$ and bound the difference by the derivative of order $j+1$
from \eqref{eq:elliptic-associated-homogeneity}. This gives
\[
 |\nabla^j v(x)|\lesssim_f |x|^{\nu-j-1}(1+\log|x|)
 \qquad(|x|\gg 1).
\]
For large $R$, the commutator $Pu_R-f=[P,\chi(\cdot/R)]v$ is supported in
$\{R\le|x|\le 2R\}$ and bounded there by $C_f R^{-d-1}(1+\log R)$, whence
$\|Pu_R-f\|_1\to 0$. On the same annulus, the commutator
$[L,\chi(\cdot/R)]v=Lu_R-\chi(\cdot/R)(K*f)$ is bounded by
$C_f R^{\nu-k-1}(1+\log R)$ (and vanishes identically when $k=0$). Since
$p(\nu-k)=-d$, this error tends to zero in $L^p$. The bound for $j=k$ also
gives $K*f\in L^p$, so $Lu_R\to K*f$ in $L^p$. Passing to the limit via
\eqref{eq:spatial-vector-Lp-continuity} gives
$C_{b,0,\mathbf 1}^*\le C_b^{P,L,*}$. Applying
Corollary~\ref{cor:weighted-homogeneous-target}\textup{(i)} and
\eqref{eq:continuous-dyadic-equivalence} completes the proof.
\end{proof}

\begin{corollary}
\label{cor:laplacian-measurable-target}
Let $K_\Delta$ denote the gradient of the radial fundamental solution of
$-\Delta$. Set
\[
 K_\Delta(x):=-\frac{x}{|\mathbb S^{d-1}|\,|x|^d}\quad(x\ne0),
 \qquad p_\Delta:=\frac{d}{d-1}.
\]
Assume $\Phi:\R^d\to\R$ is positively $p_\Delta$-homogeneous with $\Phi(0)=0$
and $L_\Phi<\infty$. For $b\in L^\infty(\R^d;\R)$ with $B=\|b\|_\infty$, write
$\mathfrak P_b^\Delta:=\mathfrak P_b[K_\Delta,\Phi]$ and
$\mathfrak P_b^{\mathrm c,\Delta}:=\mathfrak P_b^{\rm c}[K_\Delta,\Phi]$. Then
\begin{equation}\label{eq:laplacian-measurable-target}
 \left|\int_{\R^d}b(x)\Phi(\nabla u(x))\,\dd x\right|
 \le C\|\Delta u\|_{L^1(\R^d)}^{p_\Delta}
 \qquad\bigl(u\in C_c^\infty(\R^d;\R)\bigr)
\end{equation}
holds with a finite constant if and only if $\mathfrak P_b^\Delta<\infty$, and
\begin{equation}\label{eq:laplacian-quantitative-comparison}
 B+C_b^{\Delta,*}\asymp_{d,\Phi} B+\mathfrak P_b^\Delta
 \asymp_{d,\Phi} B+\mathfrak P_b^{\mathrm c,\Delta}.
\end{equation}
Taking $b=\one_E$ characterizes every measurable target $E\subset\R^d$, with
$1+C_E^{\Delta,*}\asymp_{d,\Phi}1+\mathfrak P_E^\Delta$.
\end{corollary}

\begin{proof}
Apply Corollary~\ref{cor:elliptic-measurable-target} with $P=-\Delta$ and
$L=\nabla$.
\end{proof}

\begin{remark}\label{rem:spatial-elliptic-integrands}
The same cutoff argument applies to a spatial integrand $\Psi(x,Lu)$ obeying
\eqref{eq:spatial-nonlinearity-norm}, with scalar source $f=Pu$. It yields
$L_\Psi+C^*\asymp L_\Psi+\mathfrak P_\Psi[K]$, where
\[
 \mathfrak P_\Psi[K]:=\sup_{a\in\R^d,\ \sigma=\pm 1,\ N\le M}
 \left|\sum_{n=N}^M\int_{\R^d}\Psi\bigl(x,\sigma K_n(x-a)\bigr)\,\dd x\right|.
\]
\end{remark}

\begin{remark}\label{rem:polyharmonic-chain}
For every integer $q\ge 1$, taking $P=(-\Delta)^q$ and
$L=\nabla(-\Delta)^{q-1}$ gives $\alpha=1$ and $L\Gamma_P=K_\Delta$. Their
Fourier transforms coincide off zero, so the difference is a polynomial
homogeneous of degree $1-d<0$ and must vanish. Thus the polyharmonic estimate
\begin{equation}\label{eq:polyharmonic-chain-estimate}
 \left|\int_{\R^d}b\,\Phi\bigl(\nabla(-\Delta)^{q-1}u\bigr)\,\dd x\right|
 \le C\|(-\Delta)^q u\|_{L^1(\R^d)}^{d/(d-1)}
 \qquad\bigl(u\in C_c^\infty(\R^d;\R)\bigr)
\end{equation}
is governed by the exact same criterion $\mathfrak P_b^\Delta<\infty$ as
\eqref{eq:laplacian-measurable-target}. In particular, the bounded $C^1$
domain in Proposition~\ref{prop:c1-nondini-warning} (with $\alpha=1$ and
$\Phi(z)=z_1^2-z_2^2$) violates both \eqref{eq:laplacian-measurable-target}
and \eqref{eq:polyharmonic-chain-estimate} despite cancellation on every
tangent hemisphere. More generally, $L=D^k$ yields the polyharmonic
derivative characterization whenever $0<2q-k<d$.
\end{remark}

\subsection{Spatial coefficients, vector kernels, and nonstationary shells}
\label{sec:nonstationary-shells}

Theorem~\ref{thm:path-family-classification} directly accommodates spatially
varying coefficients such as $\Psi(x,z)=\Phi(A(x)z)$ for any bounded
measurable matrix field $A$ (with $L_\Psi\le C_p L_\Phi\|A\|_\infty^p$).
Examples also include $\Psi(x,z)=z^TQ(x)z$ for $p=2$ and finite sums
$\Psi(x,z)=\sum_j b_j(x)\Phi_j(z)$, with $Q,b_j$ bounded measurable and
$\Phi_j$ satisfying the standing assumptions. Testing a sum as one integrand
retains cancellation between its terms.

In the vector-source setting ($m\ge 2$), a stronger---though simpler---sufficient
criterion controls the spherical Lipschitz norm of cumulative directional
responses. Put $Q_{a,N,M}(v):=\sum_{n=N}^M q_{n,a}(v)$ and
\begin{align*}
 \mathfrak P_{\rm dir}&:=\sup_{a,N\le M,\,|v|=1}|Q_{a,N,M}(v)|,\\
 \mathfrak P_{\rm Lip}&:=\sup_{a,N\le M}
 \Bigl(\|Q_{a,N,M}\|_{L^\infty(\mathbb S^{m-1})}
 +\operatorname{Lip}_{\mathbb S^{m-1}}Q_{a,N,M}\Bigr).
\end{align*}
Then
\begin{equation}\label{eq:directional-Lip-sufficient}
 \mathfrak P_{\rm dir}\le 2^p\mathfrak V_\Psi^H\le C_p\mathfrak P_{\rm Lip}.
\end{equation}
The first inequality follows from constant paths; for the second, writing
$F_k:=\sum_{n=N}^k q_{n,a}$ and summing by parts gives
$\sum_{n=N}^M q_{n,a}(v_n)=F_M(v_M)+\sum_{n=N}^{M-1}[F_n(v_n)-F_n(v_{n+1})]$,
bounded by $C_p\mathfrak P_{\rm Lip}D(v_N,\ldots,v_M)^p$ via
\eqref{eq:spatial-vector-lipschitz}.

If $\Psi(x,z)=\Phi(z)$ and $\Phi$ is smooth near the origin, Taylor expansion
and positive homogeneity imply that $\Phi$ is a homogeneous polynomial of
degree $p$ for integer $p$, and zero for noninteger $p$. For polynomial $\Phi$,
the $Q_{a,N,M}$ lie in a fixed finite-dimensional space of homogeneous
polynomials. The spherical supremum norm controls the Lipschitz seminorm on
this space, so $\mathfrak P_{\rm dir}<\infty$ suffices for path testing.
Homogeneous matrix kernels on the whole space give another case where fixed
directions suffice.

\begin{corollary}
\label{cor:homogeneous-vector-cancellation}
Fix $p>1$, $\beta=d/p$, and positive diagonal weights $\mathbf a$ with
$\sum_j a_j=d$.
Let $K:\R^d\setminus\{0\}\to\mathcal L(\R^m,\R^\ell)$ be homogeneous of degree
$-\beta$ under $\delta_t^{\mathbf a}$ with Lipschitz spherical restriction.
Let $\Phi:\R^\ell\to\R$ be positively $p$-homogeneous with $\Phi(0)=0$ and
$L_\Phi<\infty$. For the source class in the mean-zero estimate
\[
 \left|\int_{\R^d}\Phi(K*f)\,\dd x\right|\le C\|f\|_1^p
 \qquad\Bigl(f\in C_c^\infty(\R^d;\R^m),\ \int_{\R^d} f=0\Bigr)
\]
the nonlinear integrals are absolutely convergent. This estimate holds
uniformly if and only if
\begin{equation}\label{eq:all-direction-angular-cancellation}
 \int_{\mathbb S^{d-1}}w_{\mathbf a}(\omega)\,\Phi\bigl(K(\omega)v\bigr)\,
 \dd\sigma(\omega)=0
 \qquad\text{for every }v\in\R^m.
\end{equation}
\end{corollary}

\begin{proof}
By Lemma~\ref{lem:aniso-polar} and \eqref{eq:aniso-dini-translation}, the
sharp shells $H_n=K_n$ are admissible and every shell response $q_{n,a}(v)$
equals $\log 2$ times the integral in
\eqref{eq:all-direction-angular-cancellation}. If $\mathfrak V_\Phi^H<\infty$,
constant paths over arbitrarily long blocks force
\eqref{eq:all-direction-angular-cancellation}. If
\eqref{eq:all-direction-angular-cancellation} holds, every shell response
vanishes, so $\mathfrak V_\Phi^H=0$. The result now follows from
Theorem~\ref{thm:path-family-classification}\textup{(i)}.
\end{proof}

The same convex-shell hypotheses permit nonstationary rectangular scales
with fixed axes and varying side lengths.

\begin{proposition}
\label{prop:nonstationary-shell-classification}
Let $s_j(n)>0$ ($1\le j\le d$) be nonincreasing in $n\in\mathbb Z$ and satisfy
$s_j(n+k)\le C_s\,2^{-b_*k}s_j(n)$ for all $k\ge 0$ and some $b_*>0$. Set
$B_n:=\prod_{j=1}^d[-s_j(n),s_j(n))$ and $V_n:=|B_n|$. If bounded Borel
kernels $H_n:\R^d\to\R^\ell$ are supported in $B_n\setminus B_{n+1}$ and obey
\begin{equation}\label{eq:nonstationary-family}
 \|H_n\|_\infty\le C_0 V_n^{-1/p},
 \qquad
 \|H_n(\cdot-h)-H_n\|_1
 \le C_1 V_n^{1/p'}\vartheta\Bigl(\max_{1\le j\le d}\frac{|h_j|}{s_j(n)}\Bigr),
\end{equation}
then all conclusions of Corollary~\ref{thm:annular-family-classification} hold
for these shells, with structural constants depending only on
$(d,\ell,p,C_s,b_*,C_0,C_1,\mathcal D_\vartheta)$ and $L_\Phi$.
\end{proposition}

\begin{proof}
Set $A_n:=\overline{B_n}$ and change shell representatives on their null
boundaries before applying Corollary~\ref{thm:annular-family-classification}.
\end{proof}

\section{Sharpness and conditional cancellation}\label{sec:sharpness}

The examples test angular path variation for vector sources, spherical
Lipschitz regularity of the nonlinearity, the source conditions for automatic
integrability, boundedness of target weights, and the Dini tangent criterion.
They also separate conditional cross-scale cancellation from absolute
shell summability and test the relative contraction and Dini translation
hypotheses on the kernels.

\subsection{Necessity of angular path variation for vector sources}

In Theorem~\ref{thm:path-family-classification}, controlling path variation is
indispensable: uniform bounds on all fixed-direction block sums
($\mathfrak P_{\rm dir}<\infty$ in \eqref{eq:directional-Lip-sufficient}) do
not suffice for vector-valued sources, even when $m=\ell=2$ and all angular
profiles are $C^\infty$.

\begin{proposition}
\label{prop:directional-testing-fails}
There exist a positively $2$-homogeneous $\Phi:\mathbb C\to\R$ with $C^\infty$
spherical restriction and an annular operator family
$H_n:\R^2\to\mathcal L(\mathbb C,\mathbb C)$ satisfying
\eqref{eq:annular-family-hypotheses} with $\alpha=1$ and
$\vartheta(t)=\min(t,1)$ such that
\begin{equation}\label{eq:directional-value-primitive-bounded}
 \sup_{\substack{a\in\R^2,\ N\le M\\|v|=1}}
 \left|\sum_{n=N}^M\int_{\R^2}\Phi\bigl(H_n(x-a)v\bigr)\,\dd x\right|\le 4,
\end{equation}
every angular profile is $C^\infty$ with uniformly bounded $C^1$ norm, and
$K_H\in L^1(\R^2)$ is compactly supported, yet there exists a sequence
$f_J\in C_c^\infty(\R^2;\mathbb C)$ supported in a fixed ball with
$\int_{\R^2}f_J=0$ and $\sup_J\|f_J\|_1<\infty$ such that
\begin{equation}\label{eq:directional-counter-growth}
 \int_{\R^2}\Phi(K_H*f_J)\,\dd x\ge 2J-C.
\end{equation}
In particular, $\mathfrak V_\Phi^H=\infty$ despite
\eqref{eq:directional-value-primitive-bounded}.
\end{proposition}

\begin{proof}
\emph{Step 1: Construction and single-shell response.}
Identify $\mathbb C$ with $\R^2$ and let operators act by complex
multiplication. For $j\ge 1$, put $k_j:=2^j+1$, $a_j:=2^{-k_j^2}$, and
$\varepsilon_j:=k_j^{-1}$, and define
\[
 \varphi(\theta):=\sum_{j\ge 1}a_j\cos(k_j\theta),
 \qquad \Phi(z):=|z|^2\varphi(\arg z)\quad(z\ne 0),
 \qquad \Phi(0):=0.
\]
Then $\varphi\in C^\infty(\mathbb S^1)$, and since every $k_j$ is odd,
$\Phi(-z)=-\Phi(z)$. On the unit annulus $\{1/2\le r<1\}$, set
\[
 F_j(re^{i\theta}):=\varepsilon_j r^{-1}
 \sqrt{1+\tfrac{1}{2}\cos(k_j\theta)}\,e^{i\theta},
\]
extended by zero elsewhere. Because $\varepsilon_j k_j=1$, the angular and
radial derivatives and boundary jumps are uniformly bounded, yielding a
uniform linear $L^1$ translation modulus. Fourier orthogonality gives
\begin{equation}\label{eq:directional-counter-shell-response}
 \int_{\R^2}\Phi\bigl(F_j(x)v\bigr)\,\dd x
 =d_j|v|^2\cos(k_j\arg v),
 \qquad d_j:=\frac{\pi}{2}(\log 2)\,\varepsilon_j^2 a_j
\end{equation}
(understood as $0$ at $v=0$).

\smallskip
\emph{Step 2: Bounded directional block sums.}
Set $N_j:=\lceil d_j^{-1}\rceil$, $T_1:=0$, and $T_{j+1}:=T_j+2N_j$. Define
$H_n(x)v:=2^n F_j(2^nx)v$ for $T_j\le n<T_j+N_j$, its negative for
$T_j+N_j\le n<T_{j+1}$, and $H_n:=0$ for $n<0$. By scaling, $(H_n)$ obeys
\eqref{eq:annular-family-hypotheses}, and $K_H\in L^1(\R^2)$ is supported in
$\overline{B(0,1)}$. By \eqref{eq:directional-counter-shell-response} and
oddness of $\Phi$, $q_n(v)=\pm d_j|v|^2\cos(k_j\arg v)$ on the two halves of
the $j$th cycle. Every complete cycle cancels identically, and for $|v|=1$
every partial sum within a cycle has magnitude at most $N_jd_j<2$, yielding
\eqref{eq:directional-value-primitive-bounded}.

\smallskip
\emph{Step 3: Divergent path response.}
For each $J\ge 1$, put $u_j:=e^{i\pi/k_j}$ and define the finite path
\[
 M_n^{(J)}:=
 \begin{cases}
  1,& T_j\le n<T_j+N_j,\ j\le J,\\
  u_j,& T_j+N_j\le n<T_{j+1},\ j\le J,\\
  0,& n<0\text{ or }n\ge T_{J+1}.
 \end{cases}
\]
Both halves contribute positively, so
$\sum_n q_n(M_n^{(J)})=2\sum_{j\le J}N_jd_j\ge 2J$, while the zero-extended
variation is bounded uniformly by $2+2\pi\sum_{j\ge 1}k_j^{-1}<\infty$.
Applying Lemma~\ref{lem:path-realization} produces bounded mean-zero sources
$g_J$ supported in $\overline{B(0,2)}$ with $\sup_J\|g_J\|_1<\infty$ and
$\int\Phi(K_H*g_J)\ge 2J-C$. Mollifying each $g_J$ at a sufficiently small
scale---using $K_H\in L^1$, Young's inequality in $L^2$, and
\eqref{eq:spatial-vector-Lp-continuity}---yields smooth sources $f_J$ supported
in $\overline{B(0,3)}$ satisfying \eqref{eq:directional-counter-growth}.
\end{proof}

The obstruction lies in the unbounded spherical Lipschitz seminorm ($\ge k_j$)
of the half-cycle primitives $N_jd_j\cos(k_j\arg v)$; this nonhomogeneous
example does not conflict with Corollary~\ref{cor:homogeneous-vector-cancellation}.

\subsection{Failure for spherical Dini nonlinearities}
\label{sec:regularity-counterexample}

The point-mass criterion can fail for a Dini-continuous spherical restriction,
even when $\Phi|_{\mathbb S^{\ell-1}}$ belongs to every $C^\gamma$ with
$0<\gamma<1$.

\begin{proposition}
\label{prop:dini-nonlinearity-obstruction}
For every $d\ge 2$ and $p>1$ ($\beta=d/p$, $\alpha=d-\beta$), there exist a
critical kernel $K(x)=|x|^{-\beta}\widetilde K(x/|x|)$ with smooth
$\widetilde K:\mathbb S^{d-1}\to\R^2$ and a nonnegative, positively
$p$-homogeneous $\Phi:\R^2\to\R$ whose spherical restriction has a Dini
modulus and belongs to $C^\gamma(\mathbb S^1)$ for every $0<\gamma<1$, such
that
\begin{equation}\label{eq:dini-phi-pointmass-zero}
 \Phi\bigl(\widetilde K(\omega)\bigr)=\Phi\bigl(-\widetilde K(\omega)\bigr)=0
 \qquad(\omega\in\mathbb S^{d-1}),
\end{equation}
yet the endpoint inequality fails both on $B(0,1)$ and on $\R^d$ for mean-zero
sources in $C_c^\infty(\R^d;\R)$ (even though every individual integral
converges absolutely).
\end{proposition}

\begin{proof}
\emph{Step 1: Regularity and vanishing point-mass profiles.}
Choose smooth nonnegative $u,v$ on $\mathbb S^{d-1}$ supported in the disjoint
caps $|\omega-e_i|<1/10$ ($i=1,2$) and equal to $1$ on $|\omega-e_i|\le 1/20$,
and set $\widetilde K:=(u,v)$. Define $F(t):=t\log(e/t)$ for $0<t\le 1$ with
$F(0):=0$, and set
\begin{equation}\label{eq:dini-phi-example}
 \Phi(z):=|z|^p\Bigl(\frac{z_1}{|z|}\Bigr)^2 F\Bigl(\frac{|z_2|}{|z|}\Bigr)
 \quad(z\ne 0),\qquad \Phi(0):=0.
\end{equation}
Since $F$ is increasing and concave on $[0,1]$ with $\int_0^1 F(t)\,\dd t/t=2$,
we have $|\Phi(\theta)-\Phi(\theta')|\le 2|\theta-\theta'|+F(|\theta-\theta'|)\le 3F(|\theta-\theta'|)$
for $\theta,\theta'\in\mathbb S^1$ with $|\theta-\theta'|\le 1$, so
$\Phi|_{\mathbb S^1}$ is Dini-continuous and lies in
$\bigcap_{0<\gamma<1}C^\gamma(\mathbb S^1)$ (though not Lipschitz at $(1,0)$).
Disjointness of $\supp u$ and $\supp v$ gives
\eqref{eq:dini-phi-pointmass-zero}, so every signed point-mass functional
vanishes identically.

\smallskip
\emph{Step 2: Blow-up on $B(0,1)$ and on $\R^d$.}
For $k\ge 2$, put $r_k:=2^{-k}$, $y_k:=r_k(e_1-e_2)$, and
$w_k:=k^{-1}(\log k)^{-2}$, so that $M:=1+\sum_{k\ge 2}w_k<\infty$. Fix
$0<c<1/100$ and set $\mathcal B_k:=B(r_ke_1,cr_k)$; these balls are pairwise
disjoint and lie inside $B(0,1)$. Let $\eta\in C_c^\infty(B(0,1))$ be a
nonnegative probability mollifier, choose $0<\varepsilon_N<cr_N/10$, and put
\[
 f_N:=\eta_{\varepsilon_N}+\sum_{k=2}^N w_k\eta_{\varepsilon_N}(\cdot-y_k),
 \qquad \|f_N\|_1\le M.
\]
For $x\in\mathcal B_k$ ($2\le k\le N$), the displacements from the origin and
from $y_k$ lie within $2cr_k$ of $r_ke_1$ and $r_ke_2$, respectively; their
normalized directions lie within $4c<1/20$ of $e_1$ and $e_2$, where $u=1$ and
$v=1$. Since all other contributions are nonnegative, $Z_N:=K*f_N$ satisfies
\begin{equation}\label{eq:dini-phi-field-lower}
 (Z_N)_1(x)\ge c_1 r_k^{-\beta},\qquad
 (Z_N)_2(x)\ge c_1 w_k r_k^{-\beta}.
\end{equation}
Conversely, $|r_ke_1-y_j|^2=(r_k-r_j)^2+r_j^2\ge r_k^2/2$ for all $j$, so the
evaluation and mollifier radii leave distance at least $r_k/2$ from every
source component (and from the origin), giving
\begin{equation}\label{eq:dini-phi-field-upper}
 c_1 r_k^{-\beta}\le |Z_N(x)|\le C_1 r_k^{-\beta}
 \qquad(x\in\mathcal B_k,\ 2\le k\le N).
\end{equation}
Combining \eqref{eq:dini-phi-field-lower}--\eqref{eq:dini-phi-field-upper}
with monotonicity of $F$ gives $\Phi(Z_N(x))\ge c_2 r_k^{-d}F(c_2 w_k)$ on
$\mathcal B_k$, whence
\[
 \int_{B(0,1)}\Phi(K*f_N)\,\dd x
 \ge c_3\sum_{k=2}^N F(c_2 w_k)\longrightarrow\infty
 \qquad\Bigl(F(c_2 w_k)\asymp\frac{1}{k\log k}\Bigr).
\]

To enforce zero mean, fix a nonnegative unit-mass $\eta_*\in C_c^\infty(B(4e_1,1))$
and set $g_N:=f_N-(\int f_N)\eta_*$, so $\int g_N=0$ and $\|g_N\|_1\le 2M$.
For $x\in B(0,1)$ and $y\in B(4e_1,1)$, $(x-y)_1/|x-y|<-1/3$, outside both
angular caps; hence $K*\eta_*=0$ on $B(0,1)$, and nonnegativity of $\Phi$
gives $\int_{\R^d}\Phi(K*g_N)\ge\int_{B(0,1)}\Phi(K*f_N)\to\infty$.
For each smooth compactly supported mean-zero $g$, the field is locally
bounded and satisfies $|K*g(x)|\lesssim_g |x|^{-\beta-1}$ at infinity.
Since $|\Phi(z)|\le |z|^p$ and $\beta p=d$, each integral is finite.
\end{proof}

\subsection{Source concentration, noncompact support, and singular weights}

Local $L^p$ integrability of the potential and compact support of the source
ensure absolute integrability for zero mass. These examples locate the
possible failures at the singularity and at infinity. They also test the
boundedness assumption on target weights.

\begin{proposition}
\label{prop:l1-source-integrability-obstruction}
Even for a smooth angular kernel and a smooth target with
$\mathfrak P_E<\infty$, a compactly supported mean-zero $L^1$ source can
produce a nonlinear integrand whose positive and negative parts both have
infinite integral.
\end{proposition}

\begin{proof}
Take $d=2$, $\alpha=1$, $p=2$, $K(x)=x/|x|^2$, $\Phi(z)=z_1^2-z_2^2$, and
$E=B(0,1)$. Both signed angular profiles equal $\cos(2\theta)$, which cancels
on every semicircle, so $\mathfrak P_E<\infty$ by
Corollary~\ref{thm:smooth-domain-recovery}. Fix $0<\gamma\le 1/2$ and put
\[
 h_\gamma(y):=\frac{\gamma\,\one_{\{0<|y|<e^{-1}\}}}
 {2\pi|y|^2[\log(1/|y|)]^{1+\gamma}}\in L_c^1(\R^2),
\]
which is radial, nonnegative, and has integral $1$. Let
$h_{\rm far}\in C_c^\infty(\{2<|y|<3\})$ be radial, nonnegative, with integral
$1$, and set $f:=h_\gamma-h_{\rm far}$. For any radial $h$, circular averaging
of $K(x-y)$ gives $K*h(x)=M_h(|x|)x/|x|^2$, where
$M_h(r)=\int_{|y|<r}h(y)\,\dd y$. On $\{0<r<e^{-1}\}$,
$M_{h_\gamma}(r)=[\log(1/r)]^{-\gamma}$ and $K*h_{\rm far}=0$, whence
\[
 \Phi\bigl(K*f(r\cos\theta,r\sin\theta)\bigr)
 =\frac{\cos(2\theta)}{r^2[\log(1/r)]^{2\gamma}}.
\]
Both positive and negative parts have integral over $E$ at least
$2\int_0^{e^{-1}}r^{-1}[\log(1/r)]^{-2\gamma}\,\dd r=\infty$.
\end{proof}

Without compact support the same pair $(K,\Phi)$ gives failure at infinity.
Choose a smooth radial mass profile $M(r)$ vanishing near $0$ and equal to
$(\log r)^{-1/2}$ for $r\ge e^2$. Then $f(x):=M'(|x|)/(2\pi|x|)$ is smooth,
bounded and integrable since $\|f\|_1=\int_0^\infty|M'(r)|\,\dd r<\infty$.
Its mass is zero because $M(0)=M(\infty)=0$. The radial convolution formula
gives the bounded field $K*f(r\omega)=M(r)\omega/r$. For $r\ge e^2$ the
nonlinear density is $(\omega_1^2-\omega_2^2)/(r^2\log r)$, whose positive and
negative integrals both diverge despite vanishing whole-space shell totals.

For any finite vector measure with integrable nonlinear density,
\eqref{eq:absolute-density-decoupling} gives absolute $L^1$ summability of
the shell densities. Passing to the limit in \eqref{eq:path-block-transfer}
and adding the remainder in \eqref{eq:spatial-vector-density-series} gives
the path bound for this measure.

In the mean-zero branch, equality of optimal constants can fail on enlarging
the source class to finite measures with integrable nonlinear density.
In dimension two let
$F(x)=x/|x|$, $K=DF$ off the origin with $K(0)=0$, and $\Phi(Z)=\det Z$ on
real $2\times2$ matrices. Then $K(r\omega)=r^{-1}(I-\omega\otimes\omega)$.
The isotropic sharp shells are admissible for $p=2$ and all point-mass
responses vanish. For every finite real signed measure $\nu$, the density
is absolutely integrable and
\[
 \int_{\R^2}\det(K*\nu)\,\dd x
 =\pi\left[\nu(\R^2)^2-\sum_a\nu(\{a\})^2\right].
\]
Indeed, for distinct $y,z$ put
$C_{y,z}(x)=\det(K(x-y)+K(x-z))$. The two matrices are positive
semidefinite, so $C_{y,z}\ge0$. For $u=F(\cdot-y)+F(\cdot-z)$, apply
$\det Du\,\dd x_1\wedge\dd x_2=\dd(u_1\,\dd u_2)$ on a large disk
with small disks about $y,z$ removed. On the outer circle $u$ tends in
$C^1$ as a function of its angular parameter to $2\omega$, and on each
inner circle to $\omega$ plus a constant. The counterclockwise
boundary integrals tend to $4\pi,\pi,\pi$, respectively. Hence
$\int_{\R^2} C_{y,z}\,\dd x=2\pi$. Set $C_{y,y}=0$. The convolution is
absolutely defined almost everywhere since $K\in L^1+L^\infty$. Quadratic
polarization gives
$\det(K*\nu)=\tfrac12\iint C_{y,z}\,\dd\nu(y)\,\dd\nu(z)$ almost everywhere,
and $\int C_{y,z}\le2\pi$ gives an absolute bound of $2\pi\|\nu\|_{\rm TV}^2$
for the integral against $\dd x\,\dd|\nu|(y)\,\dd|\nu|(z)$. This justifies
Fubini.
The diagonal in the product integral contributes precisely the sum of
the squares of the atomic masses to $\nu(\R^2)^2$.

For mean-zero $\nu$, the positive and negative masses each equal
$\|\nu\|_{\rm TV}/2$, so the modulus of the integral is bounded by
$\pi\|\nu\|_{\rm TV}^2/2$, with equality at $\nu=\delta_y-\delta_z$
for $y\ne z$. The signed-integral estimate therefore has best constant
$\pi/2$ on all mean-zero measures. It has best constant zero on smooth
compactly supported mean-zero densities and hence on the mean-zero part of
$\mathcal X_H$ by
Theorem~\ref{thm:path-family-classification}. The two-atom field lies outside
$L^2_{\mathrm{loc}}$ because of its $1/|x-y|$ singularity at $y$.
Mollifying the two atoms keeps the integrals at zero, whereas the
integral of the limiting density equals $-2\pi$.
Without the mean-zero constraint the optimal constant is $\pi$ on both
classes. Indeed, $\nu(\R^2)^2$ and $\sum_a\nu(\{a\})^2$ each lie in
$[0,\|\nu\|_{\rm TV}^2]$, and a smooth nonnegative probability density
attains the bound.

We next display a nontrivial universal pair and use it to show that the
$L^\infty$ bound on target weights is essential even when all point-mass
profiles vanish identically.

\begin{proposition}
\label{prop:nontrivial-universal-pair}
For every $d\ge 2$ and $0<\alpha<d$ ($p=d/(d-\alpha)$), take $\ell=2$,
$\widetilde K^\circ(\omega):=(\omega_1,|\omega_1|)$,
$K^\circ(x):=|x|^{\alpha-d}\widetilde K^\circ(x/|x|)$, and
$\Phi^\circ(z,s):=|s|^p-|z|^p$. Then
$\Phi^\circ(\pm\widetilde K^\circ(\omega))\equiv 0$, so $(K^\circ,\Phi^\circ)$
satisfies the uniform universal-target property
(Corollary~\ref{thm:universal-target-classification}), yet there exist a
nonnegative $f\in C_c^\infty(\R^d)$ and $r_*>0$ such that
$\int_{B(0,r_*)}\Phi^\circ(K^\circ*f)\,\dd x>0$.
\end{proposition}

\begin{proof}
For a nonnegative radial unit-mass $f_0\in C_c^\infty(\{1<|y|<2\})$, oddness
of $\omega_1$ yields
\begin{equation}\label{eq:nontrivial-universal-witness-at-zero}
 (K^\circ*f_0)(0)=\biggl(0,\int_{\R^d}|y|^{\alpha-d}\frac{|y_1|}{|y|}f_0(y)\,\dd y\biggr)
 =:(0,s_0),\qquad s_0>0.
\end{equation}
Hence $\Phi^\circ((K^\circ*f_0)(0))=s_0^p>0$, and continuity gives positivity
on a small ball $B(0,r_*)$.
\end{proof}

\begin{proposition}
\label{prop:unbounded-target-density}
There exist a homogeneous kernel $K:\R^2\setminus\{0\}\to\R^2$ of degree $-1$
with Lipschitz spherical restriction and a nonnegative, positively
$2$-homogeneous $\Phi$ with $L_\Phi<\infty$ such that $\Phi(\pm K(x))=0$
everywhere, yet
\[
 \int_{B(0,1)}|x|^{-\gamma}\Phi(K*f(x))\,\dd x\le C\|f\|_1^2
 \qquad\Bigl(f\in C_c^\infty(\R^2),\ \int_{\R^2}f=0\Bigr)
\]
fails for every $0<\gamma<2$ (and likewise with target measure $\delta_0$),
even though every integral is finite and every point-mass shell profile
vanishes.
\end{proposition}

\begin{proof}
Take $K=K^\circ$ from Proposition~\ref{prop:nontrivial-universal-pair} ($d=2$,
$\alpha=1$) and $\Phi(z):=(|z_2|-|z_1|)_+^2$. Starting from $f_0$ in
\eqref{eq:nontrivial-universal-witness-at-zero}, subtract a smooth unit-mass
function supported in $B(Re_2,1)$ with $R$ large enough that the resulting
mean-zero $g\in C_c^\infty(\R^2)$ satisfies $F(x):=\Phi(K*g(x))\ge c_0>0$ on
some ball $B_\delta$. Rescaling $g_\epsilon(x):=\epsilon^{-2}g(x/\epsilon)$
preserves $\|g_\epsilon\|_1=\|g\|_1$ and gives
$\Phi(K*g_\epsilon(x))=\epsilon^{-2}F(x/\epsilon)$. For $\epsilon\delta<1$,
\[
 \int_{B(0,1)}|x|^{-\gamma}\Phi(K*g_\epsilon(x))\,\dd x
 \ge c_0\epsilon^{-\gamma}\int_{B_\delta}|y|^{-\gamma}\,\dd y\longrightarrow\infty
 \qquad(\epsilon\downarrow 0),
\]
while for $\delta_0$ the value is $\epsilon^{-2}F(0)\to\infty$.
\end{proof}

\subsection{Boundary regularity and conditional cross-scale cancellation}

We close with geometric examples illustrating the role of boundary regularity
and the strict separation between bounded signed blocks and absolute shell
summability.

\begin{remark}
\label{rem:alternating-shell-geometry}
If two spherical sections $G_0,G_1\subset\mathbb S^{d-1}$ satisfy
\[
 \int_{G_0}\Phi(\epsilon\widetilde K)\,\dd\sigma
 =-\int_{G_1}\Phi(\epsilon\widetilde K)\,\dd\sigma=c\ne 0,
\]
then placing $G_{n\bmod 2}$ on each dyadic shell $\{2^{-n-1}\le|x-a|<2^{-n}\}$
($n\ge 0$) produces alternating shell coefficients $\pm c\log 2$ at $a$, whose
consecutive block sums are uniformly bounded even though their absolute sum
diverges. This is only a single-center local model; a full all-center
construction is given in Proposition~\ref{prop:sparse-alternating-ball-target}.
\end{remark}

\begin{remark}\label{rem:polyhedral-corner-warning}
For $d=2$, $0<\alpha<2$, $p=2/(2-\alpha)$, $\widetilde K(\omega)=\omega$, and
$\Phi(z)=2z_1z_2|z|^{p-2}$ for $z\ne0$ with $\Phi(0)=0$, both signed
profiles equal $\sin(2\theta)$, which cancels on every semicircle.
At the corner $a=0$ of $\Omega=(0,1)^2$, every shell $n\ge 0$ has
$b_n^\Omega(0,\epsilon)=(\log 2)\int_0^{\pi/2}\sin(2\theta)\,\dd\theta=\log 2$,
so $\mathfrak P_\Omega=\infty$. Thus every-direction hemisphere cancellation
does not control corners.
\end{remark}

For the planar Newton kernel, \cite[Theorem~1.6]{LongZouCumulative} constructs
a convex $C^1$ counterexample for general nonzero $2$-homogeneous integrands
with Lipschitz angular part and half-circle cancellation, together with
explicit density and domain PDE tests. Below, a one-sided logarithmic boundary
graph patched into a disk yields a direct exact computation for every
$0<\alpha<2$.

\begin{proposition}
\label{prop:c1-nondini-warning}
For every $0<\alpha<2$ ($p=2/(2-\alpha)$), there exist a critical planar pair
$(K,\Phi)$ and a bounded $C^1$ domain $\Omega\subset\R^2$ such that both
full-circle and every-semicircle integrals vanish, yet
$\mathfrak P_\Omega=\infty$ and the endpoint estimate fails on $\Omega$.
\end{proposition}

\begin{proof}
Set $\widetilde K(\omega):=\omega$ and $\Phi(z):=(z_1^2-z_2^2)|z|^{p-2}$
($z\ne 0$), so both signed profiles equal $\cos(2\theta)$, which cancels on
every semicircle. Define $h_{\rm loc}(x):=x/\log(e/x)$ for $0<x<1$ and
$h_{\rm loc}(x):=0$ for $-1<x\le 0$. Then
$h_{\rm loc}'(x)=\log(e/x)^{-1}+\log(e/x)^{-2}\to 0$ as $x\downarrow 0$, so
$h_{\rm loc}\in C^1((-1,1))$ with $h_{\rm loc}'(0)=0$, though $h_{\rm loc}'$
is not Dini-continuous at $0$. Fix $0<r_1<1/4$ and
$\chi\in C_c^\infty((-2r_1,2r_1))$ with $0\le\chi\le 1$ and $\chi\equiv 1$ on
$[-r_1,r_1]$. Setting $q(x):=1-\sqrt{1-x^2}$ and
$g(x):=\chi(x)(h_{\rm loc}(x)-q(x))$ ($|x|<1$, extended by zero), the global
$C^1$ diffeomorphism $F(x,y):=(x,y+g(x))$ maps $B((0,1),1)$ onto a bounded
$C^1$ domain $\Omega$ whose lower boundary near $0$ is $y=h_{\rm loc}(x)$.

In a small ball $B_{r_*}$, $\Omega$ equals the upper half-plane minus the
strip $W:=\{x>0,\ 0<y<x/\lambda(x)\}$, where $\lambda(x):=\log(e/x)$. By
semicircle cancellation, the upper half-plane contributes zero. On $W$, the
critical density is $(x^2-y^2)/(x^2+y^2)^2$ and $y\le x/2$, so replacing
radial truncation $\{r<|z|<r_*\}$ by $\{r<x<r_*/2\}$ alters the strip integral
by $O(1)$ (the changed regions satisfy $x\asymp r$ or $x\asymp r_*$ with
vertical integrals $O(1/x)$). Evaluating the exact vertical integral gives
\[
 \int_0^{x/\lambda(x)}\frac{x^2-y^2}{(x^2+y^2)^2}\,\dd y
 =\biggl[\frac{y}{x^2+y^2}\biggr]_0^{x/\lambda(x)}
 =\frac{1}{x\lambda(x)(1+\lambda(x)^{-2})}
 =\frac{1}{x\lambda(x)}+O\biggl(\frac{1}{x\lambda(x)^3}\biggr).
\]
Integrating over $r<x<r_*/2$ shows that the signed annular integral on $\Omega$
is $-\log\log(e/r)+O(1)\to-\infty$ as $r\downarrow 0$. Hence
$\mathfrak P_\Omega=\infty$ by Lemma~\ref{lem:continuous-dyadic-cocycle}, and
Theorem~\ref{thm:sharp-measurable-characterization} yields failure of the
endpoint inequality.
\end{proof}

\begin{example}
\label{ex:alternating-radial-family}
Set $K_n^{\rm rad}(x):=\rho_{\mathbf a}(x)^{-\beta}\one_{\{2^{-n-1}\le\rho_{\mathbf a}(x)<2^{-n}\}}$
and $\Phi(t):=t|t|^{p-1}$. For $H_n:=(-1)^n K_n^{\rm rad}$,
Lemma~\ref{lem:aniso-polar} yields $\int\Phi(\sigma H_n)=\sigma(-1)^n c_d$
with $c_d:=|\mathbb S^{d-1}|\log 2$ (using $\sum_j a_j=d$). Hence
$\Lambda_H=c_d<\infty$ despite nonzero individual shell totals, so the
whole-space mean-zero estimate holds by
Corollary~\ref{thm:annular-family-classification}. Truncating to coarse shells
$G_n:=H_n\one_{\{n<0\}}$ yields $S_{\rm fine}=0$, $S_{\rm all}=\infty$, and
$\Lambda_G=c_d$: every finite-measure target admits the unconstrained estimate
(Corollary~\ref{prop:family-universal-targets}\textup{(i)}), and the
whole-space mean-zero estimate holds, yet on the infinite-measure union of
even shells $n<0$ the center-zero test diverges and the mean-zero estimate
fails.
\end{example}

\begin{example}
\label{ex:infinite-target-absolute-tail}
For $d=2$, $\alpha=1$, $p=2$, $K(x)=x/|x|^2$, $\Phi(z)=z_1^2-z_2^2$, and
$b\equiv 1$, we have $\mathfrak P_1^H=\Lambda_H=0$ but
$\mathfrak T_1^H(0)=\infty$. For any nonnegative smooth radial unit-mass $f$
supported in $B(0,1)$, harmonicity of $K$ gives $K*f=K$ on $\{|x|>1\}$, so
both positive and negative parts of $\Phi(K*f)$ have infinite integral on
$\R^2$ despite the vanishing signed point-mass seminorm.
\end{example}

\begin{remark}\label{rem:rough-universal-target}
In the universal regime (Corollary~\ref{thm:universal-target-classification}),
even $E=S\times(0,1)^{d-1}$, where $S\subset[0,1]$ is a fat Cantor
(Smith--Volterra--Cantor) set, satisfies $\mathfrak P_E=0$. Its closure
is nowhere dense and $|E|>0$. Any open set agreeing with $E$ modulo null sets
would be nonempty and would contain a ball disjoint from $\overline E$, which
is impossible.
\end{remark}

Finally, we construct a non-universal Euclidean pair and a bounded open set
whose signed point-mass seminorm is uniformly bounded across all centers and
scales while its absolute shell sum diverges. Put
\begin{equation}\label{eq:sparse-critical-density}
 \chi(\omega):=\omega_1^2-\frac{1}{d},
 \qquad h(x):=|x|^{-d}\chi(x/|x|).
\end{equation}
Because $\chi$ is a spherical harmonic of degree $2$, the Kelvin transform
identity implies that $h$ is harmonic on $\R^d\setminus\{0\}$ for every
$d\ge 2$.

\begin{proposition}
\label{prop:sparse-alternating-ball-target}
For every $d\ge 2$ and $0<\alpha<d$ ($p=d/(d-\alpha)$), there exist an explicit
pair $(K^\sharp,\Phi^\sharp)$ satisfying the standing assumptions and a
bounded open set $E^\sharp\subset\R^d$ of finite positive measure such that
$\Phi^\sharp(\pm\widetilde K^\sharp)\not\equiv 0$,
\begin{equation}\label{eq:sparse-good-target}
 \mathfrak P_{E^\sharp}[K^\sharp,\Phi^\sharp]<\infty,
 \qquad
 \sum_{n\in\mathbb Z}\bigl|b_n^{E^\sharp}(0,+1)[K^\sharp,\Phi^\sharp]\bigr|=\infty,
\end{equation}
and $E^\sharp$ is not equal modulo null sets to any bounded open set with
$C^1$ boundary.
\end{proposition}

\begin{proof}
\emph{Step 1: Construction and alternating center-zero response.}
Set $\widetilde K^\sharp(\omega):=\omega$,
$K^\sharp(x):=|x|^{\alpha-d}x/|x|$, and
$\Phi^\sharp(z):=|z|^p\chi(z/|z|)$ ($z\ne 0$), with $\Phi^\sharp(0):=0$. Both
signed point-mass densities equal $h$ in \eqref{eq:sparse-critical-density},
and evenness of $\chi$ together with $\int_{\mathbb S^{d-1}}\chi\,\dd\sigma=0$
gives $\int_H\chi\,\dd\sigma=0$ for every hemisphere $H$.

Fix $\zeta:=10^{-2}$ and for $n\ge 1$ define
\[
 u_n:=\begin{cases}e_1,&n\text{ even},\\ e_2,&n\text{ odd},\end{cases}
 \qquad R_n:=\frac{3}{4}\,2^{-n},
 \qquad \delta_n:=\zeta\bigl(n|\chi(u_n)|\bigr)^{-1/d},
 \qquad r_n^\sharp:=R_n\delta_n.
\]
Since $|\chi(u_n)|\ge 1/d$ and $d^{1/d}\le 2$, we have
$\delta_n\le 2\zeta<1/20$, so the balls $B_n:=B(R_nu_n,r_n^\sharp)$ lie
strictly inside their respective dyadic annuli $\{2^{-n-1}\le|x|<2^{-n}\}$ and
are pairwise disjoint. Put $E^\sharp:=\bigcup_{n\ge 1}B_n$; then
$|B_n|=v_d\zeta^d R_n^d/(n|\chi(u_n)|)$ ($v_d=|B(0,1)|$) and
$0<|E^\sharp|<\infty$. By harmonicity of $h$ and the ball mean-value property,
\begin{equation}\label{eq:sparse-alternating-asymptotic}
 b_n^{E^\sharp}(0,+1)[K^\sharp,\Phi^\sharp]
 =\int_{B_n}h(x)\,\dd x=|B_n|h(R_nu_n)
 =v_d\zeta^d\frac{(-1)^n}{n},
\end{equation}
so the absolute shell series diverges while all block sums at $a=0$ are
uniformly bounded.

\smallskip
\emph{Step 2: Uniform all-center control and density obstruction.}
Fix any center $a\ne 0$ with $A:=|a|$, and write
$h_n:=h\one_{\{2^{-n-1}\le|x|<2^{-n}\}}$. The bounded unit shell $h_0$
is smooth inside the annulus with jumps only across its two smooth
boundary spheres, so it belongs to $BV$. The $BV$ translation estimate from the
proof of Lemma~\ref{cor:angular-dini-sufficient} and the scaling
$h_n(x)=2^{dn}h_0(2^nx)$ give $\|h_n(\cdot-a)-h_n\|_1\lesssim 2^n A$.
\begin{itemize}
\item \emph{Coarse scales ($2^{-n}\ge 4A$):} the translation errors sum
geometrically, so \eqref{eq:sparse-alternating-asymptotic} bounds all coarse
blocks uniformly.
\item \emph{Intermediate scales ($A/4<2^{-n}<4A$):} only $O(1)$ shells occur,
contributing $O(1)$ by the uniform $L^1$ shell bound.
\item \emph{Fine scales ($2^{-n}\le A/4$):} only balls $B_j$ meeting
$B(a,A/4)$ contribute; any such $B_j$ has $R_j\in[A/2,2A]$, so at most $O(1)$
balls occur. By Corollary~\ref{thm:smooth-domain-recovery} and
Lemma~\ref{lem:continuous-dyadic-cocycle}, the unit ball has a uniformly
bounded continuous annular seminorm, which transfers to each $B_j$ by scaling
and translation invariance of the degree-$(-d)$ kernel $h$.
\end{itemize}
Combining the three ranges gives $\mathfrak P_{E^\sharp}[K^\sharp,\Phi^\sharp]<\infty$.

Finally, every neighborhood of the origin contains some $B_n$, whereas
\[
 \frac{|E^\sharp\cap B(0,r)|}{|B(0,r)|}
 \lesssim_d \frac{1}{\log(2/r)}\longrightarrow 0
 \qquad(r\downarrow 0).
\]
Any open representative agreeing with $E^\sharp$ modulo null sets must
therefore have $0$ on its boundary with density $0$. At a boundary point
where the boundary is locally a $C^1$ graph, at least one side of the graph
lies in the open set, giving density at least $1/2$.
\end{proof}

\begin{remark}
\label{rem:rough-example-scope}
Set
\[
 \mathfrak P_E^{\rm abs}[K,\Phi]
 :=\sup_{a\in\R^d,\ \epsilon=\pm 1}
 \sum_{n\in\mathbb Z}\bigl|b_n^E(a,\epsilon)[K,\Phi]\bigr|.
\]
We always have $\mathfrak P_E\le\mathfrak P_E^{\rm abs}$, whereas
Proposition~\ref{prop:sparse-alternating-ball-target} proves that finiteness
of $\mathfrak P_E$ does not imply finiteness of $\mathfrak P_E^{\rm abs}$ even
uniformly over all centers and radial truncations.
\end{remark}

\subsection{Relative contraction and the Dini assumption}

Both examples take $d=2$, $m=\ell=1$, $p>1$ and
$\Psi(x,z)=\Phi(z):=z|z|^{p-1}$. Fix a nonnegative radial probability
density $\chi\in C_c^\infty(B(0,10^{-3}))$ and use the same mean-zero source
$f=\chi-\chi(\cdot-e_1)$. Only shells $H_{-k}$ with $k\ge2$ are nonzero.
In each example the kernel is odd in one coordinate, so all point-mass
responses vanish by integration on the whole space and
$\mathfrak V_\Psi^H=0$. The coarse-only kernel is bounded, so $K_H*f$ is
locally in $L^p$.

\smallskip
\noindent\emph{Volume decay without relative contraction.}
Take $A_{-k}=[-k-1,k+1]\times[-2^k,2^k]$ for $k\ge0$ and
$A_n=[-2^{-n},2^{-n}]^2$ for $n\ge0$. The bodies are nested and
$V_{n+j}/V_n\le2^{-j}$ for $j\ge0$, since the vertical side halves at each
step. The horizontal half-width ratio of $A_0$ to $A_{-k}$ is $1/(k+1)$,
so \eqref{eq:common-scale-hypotheses} fails. Set
$H_{-k}=V_{-k}^{-1/p}\one_{E_{-k}}\operatorname{sgn}(x_2)$ for $k\ge2$.
The shell collar estimate \eqref{eq:shell-collar} uses nesting and convexity
alone. The sign change under translation contributes at most twice the
volume of $A_{-k}\cap\{|x_2|\le|h_2|\}$. Together these estimates give
\eqref{eq:common-shell-hypotheses} with $\vartheta(t)=\min(1,t)$.

For $k\ge2$ and
$x\in S_k:=(-k-3/4,-k-1/4)\times((3/4)2^k,(7/8)2^k)$, the points
$x-y$ and $x-e_1-y$, for $y\in\supp\chi$, lie in $E_{-k}$ and
$E_{-k-1}$, respectively, with positive second coordinate. Hence
$(K_H*f)(x)=V_{-k}^{-1/p}-V_{-k-1}^{-1/p}\gtrsim V_{-k}^{-1/p}$.
The integral of the positive density over $S_k$ is bounded below by
$c|S_k|/V_{-k}\gtrsim(k+1)^{-1}$. Reflection in $x_2=0$ gives the
same lower bound for the negative part on the reflected strips, so both
parts have infinite integral.

\smallskip
\noindent\emph{Continuous translations without the Dini condition.}
Now take $A_n=[-2^{-n},2^{-n}]^2$ for every $n\in\mathbb Z$. For $k\ge2$ set
\[
 R_k=2^k\bigl([-1/2,1/2]\times[3/4,3/4+1/(8k)]\bigr)\subset E_{-k},
 \qquad H_{-k}(x)=V_{-k}^{-1/p}\one_{R_k}(x)\sin(\pi x_1).
\]
For $t=D_{-k}(h)$, the sine difference on the overlap and the rectangular
collar give
\[
 V_{-k}^{-1/p'}\|H_{-k}(\cdot-h)-H_{-k}\|_1
 \lesssim k^{-1}\min(1,2^kt)+\min(k^{-1},t).
\]
For $0<t<1/2$, split at $k=\tfrac12\log_2(1/t)$. Below this value the
first term is at most $\sqrt t$, and above it the term is
$O(1/\log(e/t))$. Thus a uniform continuous modulus is
$\vartheta(t)=1/\log(e/t)$ for $0<t\le1$, with $\vartheta(0)=0$ and
$\vartheta=1$ for $t\ge1$. It is non-Dini.

Restrict to $x_1\in[-2^{k-1}+3/2,2^{k-1}-1/2]$ and the middle half of
the vertical interval of $R_k$. For every $k\ge2$ this rectangle and its
translate by $-e_1$ stay more than $10^{-3}$ from the boundary of $R_k$.
Evenness of $\chi$ and the sign reversal of the sine under the shift give
$(K_H*f)(x)=2V_{-k}^{-1/p}\sin(\pi x_1)\int\chi(y)\cos(\pi y_1)\,\dd y$,
with the last integral positive. The horizontal interval has length
$2^k-2$ and contains complete periods. Thus each sign of the density
contributes at least $c/k$ on this rectangle. The rectangles are disjoint,
so neither part is integrable. Taking $h=e_1$ on the same rectangle also
gives a lower bound $c/k$ for the normalized translation norm at
$t=2^{-k}$. The dyadic Dini sum therefore diverges for any modulus bounding
these translations.

\section*{Declaration of AI use}
The authors used the models ChatGPT 5.6 Sol and ChatGPT 6 Sol to assist with
translating and typesetting the manuscript, locating relevant references,
checking the correctness of the mathematical arguments, and simplifying proofs.
The authors provided the proof for positive diagonal dilations.
ChatGPT 5.6 Sol suggested extending the result to general expanding
matrix dilations.
Mathematical arguments, proofs, and computations in the manuscript were
verified and finalized by the authors.

\raggedbottom

\end{document}